\documentclass[final,onefignum,onetabnum,cite]{siamart171218}

\newcommand{\beq}{\begin{equation}}
\newcommand{\eeq}{\end{equation}}
\newcommand{\beqa}{\begin{eqnarray}}
\newcommand{\eeqa}{\end{eqnarray}}

\usepackage{lipsum}
\usepackage{amsfonts}
\usepackage{graphicx}
\usepackage{epstopdf}
\usepackage{algorithmic}
\usepackage{amsmath,amssymb,fixmath}
\ifpdf
  \DeclareGraphicsExtensions{.eps,.pdf,.png,.jpg}
\else
  \DeclareGraphicsExtensions{.eps}
\fi

\newsiamremark{remark}{Remark}
\newsiamremark{hypothesis}{Hypothesis}
\crefname{hypothesis}{Hypothesis}{Hypotheses}
\newsiamthm{claim}{Claim}

\headers{A Deep Learning Framework for Elasticity Imaging}{J. Yoo, A. Wahab, and J. C. Ye}

\title{A Mathematical Framework for Deep Learning in Elastic Source Imaging\thanks{Submitted to the editors on March 5, 2018.
\funding{This work was supported by the Korea Research Fellowship Program through the National Research Foundation (NRF) funded by the Ministry of Science and ICT (NRF- 2015H1D3A1062400).}}}

\author{
JaeJun Yoo\footnotemark[2]\, \footnotemark[4] 
\and 
Abdul Wahab\footnotemark[3]\, \footnotemark[4]
\and 
Jong Chul Ye\footnotemark[4]\, \footnotemark[5]
}
\usepackage{amsopn}
\DeclareMathOperator{\diag}{diag}

\ifpdf
\hypersetup{
  pdftitle={A Mathematical Framework for Deep Learning in Elastic Source Imaging},
  pdfauthor={J. Yoo, A. Wahab, and J. C. Ye}
}
\fi

\newcommand{\LE}{\mathcal{L}_{\lambda, \mu}}
\newcommand{\T}{{t}_{\rm max}}
\newcommand{\I}{\mathcal{I}}

\newcommand{\Rc}{\mathcal{R}}
\newcommand{\K}{\kappa}
\newcommand{\ds}{\displaystyle}

\newcommand{\Rd}{{\mathbb R}}

\newcommand{\Hbc}{\boldsymbol{\cal H}}

\newcommand{\Dbc}{\boldsymbol{\cal D}}
\newcommand{\hank}{\mathbb{H}}

\newcommand{\rank}{\textrm{rank}}
\newcommand{\Sigmab}{\mathbf{\Sigma}}

\newcommand{\Cb}{{\mathbf C}}

\newcommand{\Kc}{{\mathcal K}}
\newcommand{\bA}{{\mathbf{A}}}

\newcommand{\Bb}{{\mathbf B}}
\newcommand{\gb}{{\mathbf g}}

\newcommand{\Gbc}{\boldsymbol{\cal G}}

\newcommand{\Kbc}{\boldsymbol{\cal K}}
\newcommand{\Nc}{{\mathcal N}}

\newcommand{\bd}{\mathbf{d}}
\newcommand{\fb}{{\mathbf f}}
\newcommand{\bk}{\mathbf{k}}
\newcommand{\bg}{{\mathbf{g}}}

\newcommand{\bu}{\mathbf{u}}
\newcommand{\bv}{\mathbf{v}}

\newcommand{\bx}{\mathbf{x}}
\newcommand{\by}{\mathbf{y}}
\newcommand{\bz}{\mathbf{z}}

\newcommand{\bP}{\mathbf{P}}

\newcommand{\Ub}{{\mathbf U}}
\newcommand{\Vb}{{\mathbf V}}

\def\phib{{\boldsymbol{\phi}}}
\def\psib{{\boldsymbol{\psi}}}
\newcommand{\Phib}{\mathbf{\Phi}}
\newcommand{\Psib}{\mathbf{\Psi}}
\newcommand{\Ib}{{\mathbf I}}

\newcommand{\bF}{\mathbf{F}}
\newcommand{\GG}{\mathbf{G}}
\newcommand{\bGam}{\mathbf{\Gamma}}
\newcommand{\II}{\mathbf{I}}
\newcommand{\bLambda}{{\mathbf{\Lambda}}}
\newcommand{\bS}{\mathbf{S}}

\newcommand{\NN}{\mathbb{N}}
\newcommand{\RR}{\mathbb{R}}

\begin{document}

\maketitle

\renewcommand{\thefootnote}{\fnsymbol{footnote}}

\footnotetext[2]{Clova AI Research, NAVER Corporation, Naver Green Factory, 6 Buljeong-ro, Bundang-gu, 13561, South Korea (\email{jaejun.yoo@navercorp.com}, \email{jaejun2004@kaist.ac.kr}).}
\footnotetext[3]{Department of Mathematics, University of Education, Attock Campus 43600, Attock,  Pakistan (\email{wahab@kaist.ac.kr}).}
\footnotetext[4]{Bio-imaging and Signal Processing Laboratory, Department of Bio and Brain Engineering, Korea Advanced Institute of Science and Technology, 291 Daehak-ro, Yuseong-gu, 34141, Daejeon, South Korea (\email{jong.ye@kaist.ac.kr}).}
\footnotetext[5]{Department of Mathematical Sciences, Korea Advanced Institute of Science and Technology, 291 Daehak-ro, Yuseong-gu, 34141, Daejeon, South Korea.}

\renewcommand{\thefootnote}{\arabic{footnote}}

\begin{abstract}
  An inverse elastic source problem with sparse measurements is of concern. A generic mathematical framework is proposed which extends a low-dimensional manifold regularization  in the conventional source reconstruction algorithms thereby enhancing their performance with sparse data-sets. It is rigorously established that the proposed framework is equivalent to the so-called \emph{deep convolutional framelet expansion} in machine learning literature for inverse problems.  Apposite numerical examples are furnished to substantiate the efficacy of the proposed framework.
\end{abstract}



\section{Introduction}
An abundance of real-world inverse problems, for instance in  biomedical imaging, non-destructive testing, geological exploration, and sensing of seismic events, is concerned with the spatial and/or temporal {support} localization of \emph{sources} generating wave fields in acoustic, electromagnetic, or elastic media (see, e.g., \cite{FiniteTR2011, Archer2010Passive,  Gennisson2003ASAJ, Kedar2011CRGeo,   SumatraTR2006, Lim17OE, Michel2004EEG, Porter1982JOSA,   Sabra2007passive} and references therein).   Numerous application-specific algorithms have been proposed in the recent past to procure solutions of diverse \emph{inverse source problems} from time-series or time-harmonic measurements of the generated waves (see, e.g., \cite{Albanese2006, PAT, MSRI, TRFink, Lakhal2008,   Handbook, ENoise, TREM, HongyuEM18,  JSparsity, HongyuAc18}). 
The \emph{inverse elastic source problems} are of  particular interest in this paper due to their relevence in \emph{elastography} \cite{Princetonl, Archer2010Passive,   Gennisson2003ASAJ, Sabra2007passive}. Another potential application is the localization of the background noise source distribution of earth, which contains significant information about the regional geology, time-dependent crustal changes and earthquakes \cite{Griffa2008, Kedar2011CRGeo, FiniteTRLimitations, SumatraTR2006}.

Most of the conventional algorithms are suited to continuous measurements, in other words, to experimental setups allowing to measure wave fields at each point inside a region of interest or on a substantial part of its boundary. In practice, this requires mechanical systems that furnish discrete data sampled on a very fine grid confirming to the Nyquist sampling rate. Unfortunately, this is not practically feasible due to mechanical, computational, and financial constraints. In this article, we are therefore interested in the problem of elastic source imaging with very sparse data, both in space and  time, for which the image resolution furnished by the conventional algorithms degenerates. In order to explain the idea of the proposed framework, we will restrict ourselves to the \emph{time-reversal technique} for elastic source localization presented by Ammari et al. \cite{TrElastic} as the base conventional algorithm due to its robustness and efficiency. It is precised that any other contemporary algorithm can be adopted accordingly. The interested readers are referred to the articles \cite{ TrAcoustic, TrElastic, FiniteTR2011, TRFink,FiniteTRLimitations, TREM} and reference cited therein for further details on time-reversal techniques for inverse source problems and their mathematical analysis.

%

One potential remedy to overcome the limitation of the conventional algorithms  is to incorporate the smoothness penalty such as  the total variation (TV) or other sparsity-inducing penalties under a data fidelity  term. These approaches are, however, computationally expensive due to the repeated applications of the forward solvers and  reconstruction steps during iterative updates. Direct image domain processing using these penalties could bypass the iterative applications of the forward and inverse steps, but the performance improvements are not remarkable.

 Since the deep convolutional neural network (CNN) known as AlexNet \cite{krizhevsky2012imagenet} pushed the state of the art by about 10\%, winning a top-5 test error rate of 15.3\% in the ImageNet Large Scale Visual Recognition Challenge (ILSVRC) 2012 \cite{ILSVRC15} compared to the second-best entry of 26.2\%, the performance of  CNNs  continuously improved and eventually surpassed the human-level-performance (5.1\%, \cite{ILSVRC15}) in the image classification task.
 Recently, deep learning approaches have achieved tremendous success not only for classification tasks,
 but also  in various inverse problems of computer vision area such as segmentation \cite{ronneberger2015u}, image captioning \cite{karpathy2015deep}, denoising \cite{zhang2017beyond}, and super resolution \cite{bae2017beyond},  for example. 

 {Along with those developments, by applying the deep learning techniques, a lot of studies in medical imaging area have also shown good performance in various applications \cite{adler2018learned, chen2018learn, chen2017low,  chen2017lowBOE, han2016deep, jin2016deep, kang2018deep, kang2016deep, kang2017wavresnet,       wang2016perspective, wolterink2017generative, wurfl2016deep}. 
 For example, Kang et al. \cite{kang2016deep} first successfully demonstrated wavelet domain
 deep convolutional neural network (DCNN) for low-dose computed tomography (CT), winning the second place in 2016 American Association of Physicists in Medicine (AAPM) X-ray CT Low-dose Grand Challenge \cite{mccollough2017low}. Jin et al. \cite{jin2016deep} and Han et al. \cite{han2016deep} independently showed that the global streaking artifacts from the sparse-view CT can be removed efficiently with the deep network. In MRI, Wang et al. \cite{wang2016accelerating} applied deep learning to provide a soft initialization for compressed sensing MRI (CS-MRI). In photo-acoustic tomography, Antholzer et al. \cite{antholzer2017deep} proposed a U-Net architecture  \cite{ronneberger2015u} to effectively remove streaking artifacts from inverse spherical Radon transform based reconstructed images.} 
The power of machine learning for inverse problems has been also demonstrated in material discovery and designs, in which the goal
is to find the material compositions and structures to satisfy the design goals under assorted design constraints  \cite{liu2017materials, liu2017onset, shi2015multi}.

 In spite of such intriguing performance improvement by deep learning approaches, the origin of the success  for inverse problems was poorly understood. To address this,   we recently proposed 
so-called  \emph{deep convolutional framelets}  as a powerful mathematical framework to understand deep learning approaches for inverse problems \cite{ye2017deep}. The novelty of the deep convolutional framelets was the discovery that  an encoder-decoder network structure emerges as the signal space manifestation from 
 Hankel matrix decomposition in the higher dimensional space \cite{ye2017deep}.
In addition, by controlling the number of filter channels, a neural network is trained to learn the optimal \emph{local bases} so that it gives the best low-rank shrinkage \cite{ye2017deep}.
This discovery demonstrates an important link between the deep learning and the compressed sensing approachs \cite{donoho2006compressed} through a Hankel structure matrix decomposition \cite{  jin2016general,  jin2015annihilating, ye2016compressive}.

Thus, the aim of this paper is to provide a deep learning reconstruction formula for elastic source imaging from sparse measurements. Specifically,  a generic framework is provided
that incorporates a low-dimensional manifold regularization in the conventional reconstruction frameworks. 
As it will be explained later on, the resulting algorithm can be extended to the deep convolutional framelet expansion   in order to achieve an image resolution comparable to that furnished by the continuous/dense measurements \cite{ye2017deep}.

The paper is organized as follows. The inverse elastic source problem, both in discrete and continuous settings,  is introduced in \cref{Sect:Prob} and  a brief review of the time-reversal algorithm is also provided. The mathematical foundations of the proposed deep learning approach are furnished in \cref{Sect:DL}. \Cref{Sect:NT} is dedicated to the design and training of the deep neural network. A few numerical examples are furnished in \cref{Sect:NE}. The article ends with a brief summary in \cref{Sect:Con}.

\section{Problem formulation}
\label{Sect:Prob}

Let us first mathematically formulate  the inverse elastic source problems with continuous and discrete measurements. Then, we will briefly review the time-reversal technique for elastic source imaging (with continuous data) as discussed in \cite{TrElastic} in order to make the paper self-contained.

\subsection{Inverse elastic source problem with continuous measurements}
\label{sect:model} 
Let $\bS:\RR^d\times\RR\to \RR^d$ be a compactly supported function. Then,  the wave propagation in a linear isotropic elastic medium loaded in $\RR^d$ ($d=2,3$)  is governed by the Lam\'e system,
\begin{equation*}
\left\{
\begin{array}{ll}
    \ds\frac{\partial^2 \bu}{\partial t^2 }(\bx,t) - \LE \bu(\bx,t)   =  \ds \bS(\bx,t), & (\bx,t) \in \RR^d \times \RR ,
\\
 \ds  \bu(\bx,t) = {\bf 0}= \frac{\partial \bu}{\partial t} (\bx,t), &\,\bx \in \RR^d, \,\, t
  < 0,
  \end{array}
\right.
 \end{equation*}
where $\bu=(u_1,\cdots, u_d)^\top:\RR^d\times\RR\to \RR^d$ is the elastic wave field generated by the source $\bS$, operator 
$
\LE\bu= \mu\Delta\bu+(\lambda+\mu)\nabla (\nabla\cdot\bu ) 
$
is the linear isotropic elasticity operator with Lam\'e parameters of the medium $(\lambda, \mu)$, and superscript $\top$ indicates the transpose operation. Here, it is assumed for simplicity that the volume density of the medium is unit, i.e.,  $\lambda$, $\mu$, and   $\bS$ are density normalized.   Moreover, the source is punctual in time, i.e., $\bS(\bx,t)=\bF(\bx){d\delta_0(t)}/{dt}$, where $\delta_0$ denotes the Dirac mass at $0$ and its derivative is defined in the sense of distributions.  

Let  $\Omega\subset\RR^d$ be an open bounded smooth imaging domain with $\mathcal{C}^2-$boundary $\partial\Omega$, compactly containing the spatial support  of $\bF(\bx)=(F_1,\cdots,F_d)^\top\in\RR^d$, denoted by ${\rm supp}\{\bF\}$,  i.e., there exists a compact set $\Omega^*\subset\RR^d$ strictly contained in $\Omega$ such that ${\rm supp}\{\bF\}\subset\Omega^*\subset\Omega$. 
Then, the inverse elastic source problem with continuous measurement data is to recover $\bF$ 
given the measurements  
\begin{equation*}
\Big\{ \bd(\by,t) :=  \bu(\by,t)\,\Big| \quad \forall\, \by\in\partial \Omega,\quad \forall\, t\in (0,\T) \Big\},
\end{equation*}
where $\T$ is the final control time such that $\bu(\bx,\T)\approx 0$ and $\partial_t\bu(\bx,\T)\approx 0$ for all $\bx\in\partial\Omega$. 

\begin{figure}[!htb]
\centering
\includegraphics[width=0.7\linewidth]{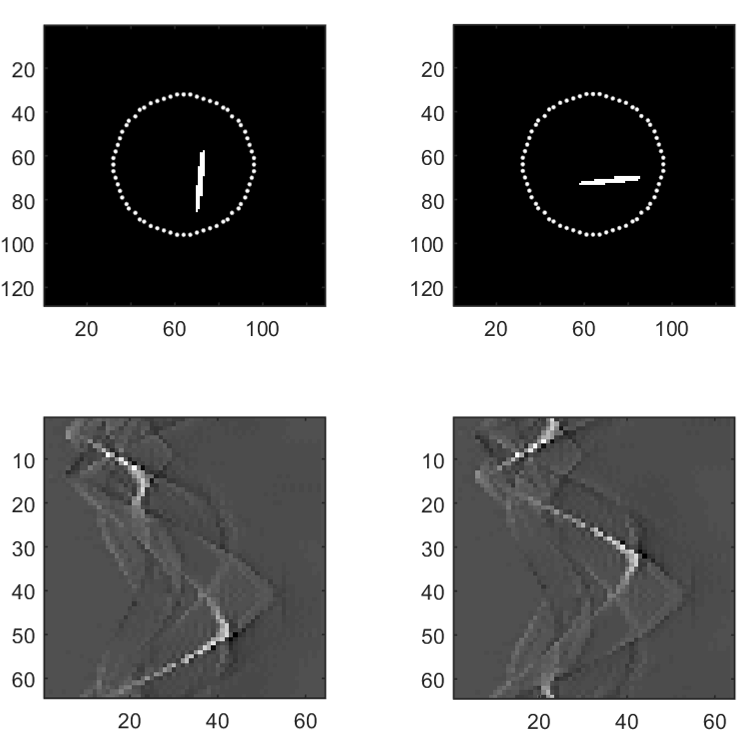}
\caption{Source and measurement configurations in 2D when the propagation direction is along the $x$-axis, the region of interest $\Omega$ is the unit disc centered at origon, $64$ detectors are placed at the control geometry $\partial\Omega$ with a time interval $[0,2s]$, the temporal scanning rate is $2^{-6}s$ and the displayed region is $[-2cm,2cm]^2$ discretized with a mesh size $2^{-7}cm$. Top: The pressure component (or $x$-component) (left) and the  shear component (or $y$-component) (right) of the spatial support of the source density $\bF$. Bottom: Measurements of the $x$-component (left) and $y$-component (right) of the wave field $\bu$ at $\partial\Omega$ (scanning times versus detector positions).}\label{fig:config}
\end{figure}

 It is precised that $\bF$ and $\bu$ can be decomposed in terms of irrotational components (or pressure components polarizing along the direction of propagation) and solenoidal components (or shear components polarizing orthogonal to the direction of propagation). In particular, in a two-dimensional (2D) frame-of-reference wherein $x$- and $y$-axes are aligned with and orthogonal to the direction of propagation, respectively, the respective components of $\bF$ are its pressure and shear components (see, e.g., \Cref{fig:config} for the imaging setup and source configuration). 

\subsection{Inverse elastic source problem with discrete measurements}

Most of the conventional algorithms require the measurement domain $\partial\Omega\times (0,\T)$ to be sampled at the Nyquist  rate so that a numerical reconstruction of the spatial support is achieved at a high resolution. Specifically,  the distance between  consecutive receivers is taken to be less than half of the wavelength corresponding to the smallest frequency in the bandwidth and the temporal scanning is done at a fine rate so that the relative difference between consecutive scanning times is very small. 

In practice, it is not feasible to place a large number of receivers at the boundary of the imaging domain and most often the measurements are available only at a few detectors (relative to the number of those required at the Nyquist sampling rate). As a result, one can not expect well-resolved reconstructed images from the conventional algorithms requiring continuous or dense measurements. 

In the rest of this subsection,  the mathematical formulation of the discrete inverse elastic source problem is provided. Towards this end, some notation is fixed upfront.  For any sufficiently smooth function $v:\RR\to\RR$, its temporal Fourier transform  is defined by
\begin{align*}
\hat{v}(\omega)=\mathcal{F}_t[v](\omega):=\int_{\RR} e^{\iota\omega t}v(t)dt,
\end{align*}
where $\omega\in\RR$ is the temporal frequency. Similarly, the spatial Fourier transform of an arbitrary smooth function $w:\RR^d\to\RR$ is defined by 
\begin{align*}
\hat{w}(\bk)=\mathcal{F}_\bx[w](\bk):=\int_{\RR^d} e^{-\iota\bk\cdot \bx}w(\bx)d\bx,
\end{align*}
with spatial frequency $\bk\in\RR^d$.
Let the function $\widehat{\GG}_{\omega}$ be the Kupradze matrix of fundamental solutions associated to the time-harmonic elastic wave equation, i.e., 
\begin{equation}
\label{eq:green0}
\LE[\widehat{\GG}_{\omega}](\bx)+\omega^2\widehat{\GG}_{\omega}(\bx)=-\delta_{\bf 0}(\bx) \II_d,\qquad \bx\in\RR^d,
\end{equation}
where $\II_d\in\RR^{d\times d}$ is the identity matrix. For later use, we decompose $\widehat{\GG}$ into its shear and pressure parts as 
\begin{equation*}
\begin{aligned}
&\widehat{\GG}_{\omega}(\bx) = \widehat{\GG}^P_{\omega}(\bx) + \widehat{\GG}^S_{\omega}(\bx), \quad\bx\neq \mathbf{0},
\\
&\widehat{\GG}^P_{\omega}(\bx) = -\frac{1}{\omega^2} \nabla\nabla^\top \widehat{g}^P_{\omega}(\bx)
\quad\text{and}\quad 
\widehat{\GG}^S_{\omega} (\bx)= \frac{1}{\omega^2} \left(
\kappa_S^2 \II_d +  \nabla\nabla^\top \right) \widehat{g}^S_{\omega}(\bx),
\end{aligned}
\end{equation*}
where  
\begin{align*}
\widehat{g}_{\omega}^\alpha(\bx)=
\begin{cases}
\ds \frac{\iota}{4}H^{(1)}_0(\K_\alpha |\bx|), & d=2,
\\
\ds\frac{1}{4\pi |\bx|}e^{i\K_\alpha |\bx|}, & d=3,
\end{cases}
\quad\text{and}\quad 
\K_\alpha:=\frac{\omega}{c_\alpha} 
\quad\text{with}\quad
\alpha=P,S.
\end{align*}
Here, $H^{(1)}_0$ denotes the  first-kind Hankel function of order zero, and  $c_P = \sqrt{\lambda + 2\mu}$ and  $c_S = \sqrt{\mu}$ are the pressure and shear wave speeds, respectively.

If $\GG(\bx,t):=\mathcal{F}^{-1}_t[\widehat{\GG}_\omega(\bx)]$  then, by invoking the Green's theorem,
\begin{align}
\bd(\by,t)=\bu(\by,t)\big|_{\by\in\partial\Omega}=\left[\int_\Omega\frac{\partial}{\partial t}\GG(\by-\bz,t)\bF(\bz)d\bz\right]\bigg|_{\by\in\partial\Omega}=:\mathcal{D}[\bF](\by,t),\label{ForwOp}
\end{align}
for all $(\by,t)\in\partial\Omega\times[0,\T]$. Here,  $ \mathcal{D}:L^2(\Omega)^d\to L^2(\partial \Omega\times  [0,T])^d$ denotes the source-to-measurement operator.

Let $\by_1,\cdots,\by_M\in\partial\Omega$  be the locations of $M\in\NN$ point receivers measuring the time-series of the  outgoing elastic wave  $\bu$ at instances $0<t_1<\cdots<t_N<\T$  for some $N\in\NN$.  Then, the inverse elastic source problem with discrete data is to recover 
$\bF$  given the discrete measurement set 
\begin{equation*}
\Big\{ \bd(\by_m,t_n) := \mathcal{D}[\bF](\by_m,t_n) \,\Big|\quad \forall\,  1\leq m\leq M,\quad 1\leq n\leq N \Big\}.
\end{equation*}
In this article, we are interested in the discrete inverse source problem with sparse data, i.e., when $M$ and $N$ are small relative to the Nyquist sampling rate. 

In order to facilitate the ensuing discussion, let us introduce the discrete measurement vector  $\Gbc\in\RR^{dMN}$ by 
\begin{align}
\Gbc:=
\begin{pmatrix}
\Gbc_1 \\ \vdots \\  \Gbc_d
\end{pmatrix},
\quad\text{where}\quad
\Gbc_i:=
\begin{pmatrix}
\Gbc^{1}_i \\ \vdots \\ \Gbc^{N}_i
\end{pmatrix}
\quad\text{with}\quad 
\Gbc^{n}_i := 
\begin{pmatrix}
[\bd(\by_1,t_n)]_i \\ \vdots \\ [\bd(\by_M,t_n)]_i
\end{pmatrix}.
\label{Gbc}
\end{align}
Here and throughout this investigation, notation $[\cdot]_i$ indicates the $i$-th component of a vector and $[\cdot]_{ij}$ indicates the $ij$-th component of a matrix.  Thus, $\Gbc^{n}_i$, for $1\leq i\leq d$, denotes the vector formed by the $i$-th components  of the waves recorded at points $\by_1,\cdots,\by_M$ at a fixed time instance $t_n$. 

Let us also introduce the forward operator, $\Dbc_{\rm dis}:L^2(\RR^d)^d\to \RR^{dMN}$, in the discrete measurement case by 
\begin{align*}
\Dbc_{\rm dis} [\bF]:= 
\begin{pmatrix}
\Dbc_1 [\bF]
\\
\vdots
\\
\Dbc _d[\bF]
\end{pmatrix},
\,\text{where}\,
\Dbc_i [\bF]= 
\begin{pmatrix}
\Dbc^{1}_i [\bF]
\\
\vdots
\\
\Dbc^{N}_i [\bF]
\end{pmatrix}
\,\text{with }\,
\Dbc^{n}_i [\bF]=
\begin{pmatrix}
\left[\mathcal{D}[\bF](\by_1,t_n)\right]_i\\\vdots\\\left[\mathcal{D}[\bF](\by_M,t_n)\right]_i
\end{pmatrix}.
\end{align*}
Then, the  inverse elastic source problem with discrete data is to recover $\bF$ from the relationship 
\begin{align}
\Gbc =\Dbc_{\rm dis}[\bF].\label{Gbc2}
\end{align}

\subsection{Time-reversal for elastic source imaging: A review}
\label{sect:TR} 

The idea of the time-reversal algorithm is based on a very simple observation that the wave operator in loss-less (non-attenuating) media is self-adjoint and that the  corresponding Green's function possesses the reciprocity property \cite{TRFink}. In other words, the wave operator is invariant under time transformation $t\to -t$ and the positions of the sources and receivers can be swapped. Therefore, it is possible to theoretically \emph{revert} a wave from the recording positions and different control times to the source locations and the initial time in chronology thereby converging to the source density. Practically, this is done by \emph{back-propagating} the measured data, after transformation $t\to \T-t$, through the adjoint waves $\bv_\tau$ (for each time instance $t=\tau$) and adding the contributions $\bv_\tau$ for all $\tau\in(0,\T)$  after evaluating them at the final time $t=\T$. Precisely, the adjoint wave $\bv_\tau$, for each $\tau\in (0,\T)$, is constructed as the solution to 
\begin{equation*}
\left\{
\begin{array}{ll}
   \ds  \frac{\partial^2 \bv_\tau}{\partial t^2 }(\bx,t) - \LE \bv_\tau(\bx,t) = \frac{d \delta_\tau(t)}{d t}\bd(\bx,\T-\tau)\delta_{\partial \Omega}(\bx), & (\bx,t)\in\RR^d\times\RR,
\\ 
 \ds  \bv_\tau(\bx,t) =  \frac{\partial\bv_\tau}{\partial t}(\bx,t)
={\bf 0}, &  \bx\in\RR^d, \, t<\tau,
\end{array}
\right.
\end{equation*}
where $\delta_{\partial \Omega}$ is the surface Dirac mass on
$\partial \Omega$. Then, the time-reversal imaging function is defined by
\begin{equation}
\label{def:I1}
 \I_{\rm TR}(\bx) = \int_0^{\T} \bv_\tau(\bx,\T) d\tau,\quad \bx\in\Omega.
\end{equation}
By the definition of the adjoint field $\bv_\tau$ and the Green's theorem,  
\begin{align*}
\bv_\tau(\bx,t)=\int_{\partial\Omega}\frac{\partial}{\partial t}\GG(\bx-\by,t-\tau)\bd(\by,\T-\tau)d\sigma(\by).
\end{align*}
Therefore, the time-reversal function can be explicitly expressed as 
\begin{align*}
\I_{\rm TR}(\bx)=\int_{0}^{\T} & \int_{\partial \Omega}\int_{\Omega}\left[\frac{\partial}{\partial t}\GG(\bx-\by,t-\tau)\right]\bigg|_{t=\T}
\nonumber
\\
&\times \left[\frac{\partial}{\partial t}\GG(\by-\bz,t)\bF(\bz)\right]\bigg|_{t=\T-\tau}d\bz d\sigma(\by) d\tau.
\end{align*} 

The time-reversal function $\I_{\rm TR}$ in \eqref{def:I1} is usually adopted to reconstruct the source distribution in an elastic medium. However,  it does not provide a good reconstruction due to a non-linear coupling between the shear and pressure parts of the elastic field $\bu$ at the boundary, especially when the sources are extended  \cite{FiniteTRLimitations, FiniteTR2011, TrElastic}. In fact, these components propagate at different wave-speeds and polarization directions, and cannot be separated at the surface of the imaging domain. If we simply back-propagate the measured data then the time-reversal operation mixes the components of the recovered support of the density $\bF$. Specifcally, it has been established in \cite{TrElastic} that, by time reversing and back-propagating the elastic wave field signals as in (\ref{def:I1}), only a blurry image can be reconstructed together with an additive term introducing the coupling artifacts.

As a simple remedy for the coupling artifacts, a surgical procedure is proposed  in \cite{TrElastic} taking the leverage of a \emph{Helmholtz decomposition} of  $\I_{\rm TR}$, (regarded as an \emph{initial guess}).   A weighted time-reversal imaging function (denoted by  $\I_{\rm WTR}$ hereinafter) is constructed by separating the  shear and  pressure components of $\I_{\rm TR}$ as
\begin{equation*}
\I_{\rm TR} =  \nabla \times \psi_{\I_{\rm TR}} + \nabla \phi_{\I_{\rm TR}},
\end{equation*}  
 and then taking their weighted sum wherein the weights are respective wave speeds and the functions $\psi_{\I_{\rm TR}}$ and $\phi_{\I_{\rm TR}}$ are obtained by solving a weak Neumann problem. Precisely,  $\I_{\rm WTR}$ is defined by 
\begin{equation}
\label{tildeI}
\I_{\rm WTR} =  c_S \nabla \times \psi_{\I_{\rm TR}} +  c_P \nabla
\phi_{\I_{\rm TR}}.
\end{equation}
In fact, thanks to the Parseval's theorem and the fact that $\bF$ is compactly supported inside $\Omega\subset\RR^d$, it can be established that
\begin{align*}
\I_{\rm WTR}(\bx)= \frac{1}{4\pi}\int_{\RR^d}\int_{\RR} \omega^2
\bigg[
&\int_{\partial \Omega} 
\bigg(
\widehat{\bGam}_\omega(\bx-\by)\overline{\widehat{\GG}_\omega(\by-\bz)}
\nonumber
\\
&+
 \overline{\widehat{\bGam}_\omega(\bx-\by)}\widehat{\GG}_\omega(\by-\bz)
\bigg)d\sigma(\by)
\bigg]d\omega \bF(\bz)d\bz,
\end{align*}
for a large final control time $\T$ with
$$
\widehat{\bGam}_\omega (\bx):=c_P\widehat{\GG}^P_\omega(\bx)+c_S\widehat{\GG}^S_\omega(\bx), \qquad \forall \,\bx\in\RR^d.
$$
After tedious manipulations, using the \emph{elastic Helmholtz-Kirchhoff identities} (see, e.g., \cite[Proposition 2.5]{TrElastic}), and assuming $\Omega$ to be a ball with radius $R\to+\infty$, one finds out that
\begin{align*}
\I_{\rm WTR} (\bx) \quad {R\to +\infty \atop =} \quad \frac{1}{2\pi}\int_{\RR^d}\int_{\RR} \omega\Im \left[\widehat{\GG}_\omega(\bx-\bz)\right] d\omega \bF(\bz)d\bz.
\end{align*}
Since 
\begin{align*}
\frac{1}{2\pi}\int_{\RR} -i\omega\widehat{\GG}_\omega(\bx-\bz) d\omega = \delta_{\bx}(\bz)\mathbf{I}_d,
\end{align*}
which comes from the integration of the time-dependent version of Eq. \eqref{eq:green0} between $t =0^-$ and $t = 0^+$, the following result holds (see, e.g., \cite[Theorem 2.6]{TrElastic}).
\begin{theorem}\label{theorem}
Let $\Omega$ be a ball in $\RR^d$ with  large radius R. Let $\bx\in\Omega$ be sufficiently far from the boundary $\partial\Omega$ with respect to the wavelength and $\I_{\rm WTR}$ be defined by \eqref{tildeI}. Then,
\begin{equation*}
\I_{\rm WTR}(\bx) \quad {R\to +\infty \atop =} \quad\bF(\bx).
\end{equation*}
\end{theorem}

We conclude this section with the following remarks. Let $\mathcal{D}$ be the source-to-measurement operator, defined in \eqref{ForwOp}. Then,  it is easy to infer  from  \cref{theorem} that its inverse (or the measurement-to-source) operator is given by 
\begin{align*}
\mathcal{D}^{-1}[\bd](\bx)  \quad {R\to +\infty \atop =} \quad \I_{\rm WTR} (\bx), 
\end{align*}
when imaging domain $\Omega$ is a ball with large radius $R$. However, there are a few technical limitations. Firstly, if  $\Omega$ is not \emph{sufficiently large} as compared to the characteristic size of the support of $\bF$, which in turn should be sufficiently localized at the \emph{center} of the imaging domain (i.e., located far away from the boundary $\partial\Omega$),  one can only get an approximation of $\bF$ which may not be very well-resolved. Moreover, $\I_{\rm WTR}$ may not be able to effectively rectify the coupling artifacts in that case   as it has been observed for extended sources in \cite{TrElastic}. Secondly, like most of the contemporary conventional techniques, time-reversal algorithm requires continuous measurements (or dense measurements at the Nyquist sampling rate). Therefore, as will be highlighted later on in the subsequent sections,  very strong streaking artifacts appear when the time-reversal algorithm is applied with sparse measurements. In order to overcome these issues,  a deep learning approach is discussed in the next section.

\section{Deep learning approach for inverse elastic source problem}\label{Sect:DL}

Let us consider the inverse elastic source problem with sparse measurements. Our aim is to recover $\bF$ from the relationship \eqref{Gbc2}. Unfortunately,   \eqref{Gbc2} is not uniquely solvable due to sub-sampling. In fact, the null space,  $\Nc({\Dbc_{\rm dis}})$, of the forward operator $\Dbc_{\rm dis}$ is non-empty, i.e., there exist non-zero functions, $\bF^0\in L^2(\RR^d)^d$, such that 
\begin{align*}
\Dbc_{\rm dis}(\bF^0)=\mathbf{0}.
\end{align*}
Moreover, the existence of the non-radiating parts of the source also makes the solution non-unique. This suggests that there are infinite many feasible solutions to the discrete problem \eqref{Gbc2}. Hence, the application of the time-reversal  algorithm requiring the availability of continuous or dense measurements results in strong imaging artifacts severely affecting the resolution of the reconstruction. 

A typical way to avoid the non-uniqueness of the solution from sparse measurements is the use of regularization. Accordingly, many regularization techniques have been proposed over the past few decades. Among various penalties for regularization,  here our discussion begins with
a low-dimensional manifold constraint using a structured low-rank penalty \cite{ye2016compressive}, which is closely related to the deep learning approach proposed in this investigation.

\subsection{Generic inversion formula under structured low-rank constraint}

Let  $\{\bz_q\}_{q=1}^Q\subset\Omega$, for some integer $Q\in\NN$, be a collection of finite number of sampling points of the region of interest $\Omega$ confirming to the Nyquist sampling rate. In this section, a (discrete) approximation of the density $\bF$ is sought using piece-wise constants or splines ansatz
\begin{align*}
[\bF(\bz)]_i:=\sum_{q=1}^Q [\bF(\bz_q)]_i \vartheta_i(\bz,\bz_q) , \quad \forall\, \bz\in\Omega,
\end{align*}
where $\vartheta_i(\cdot,\bz_q)$ is the basis function for the $i$-th coordinate, associated with $\bz_q$. Accordingly, the discretized source density to be sought is introduced by 
\begin{align*}
\fb:= 
\Big(
\fb_1^\top,
\cdots,
\fb_d^\top
\Big)^\top
\in\RR^{dQ}
\quad\text{with}\quad
\fb_i:=
\Big(
[\bF(\bz_1)]_i ,
\cdots,
[\bF(\bz_Q)]_i
\Big)^\top\in\RR^Q.
\end{align*}

{
Let us define the row-vector $\bLambda_{i,j}^{n,m}\in\RR^{1\times Q}$ by 
\begin{align*}
\left[\bLambda_{i,j}^{n,m}\right]_q:=\int_{\Omega} \left[\frac{\partial}{\partial t}\GG(\by_m-\bz,t)\right]_{ij}\bigg|_{t=t_n}\vartheta_j(\bz,\bz_q)d\bz,
\end{align*}
where superposed $n$ and $m$ indicate the dependence on $n$-th time instance for $1\leq n\leq N$ and  $m$-th boundary point $\by_m$ for  $1\leq m \leq M$, respectively. The subscripts  $1\leq i,j\leq d$ indicate that the $(i,j)$-th component of the Kupradze matrix is invoked and the  index $1\leq q\leq Q$ indicates that the basis function associated with the internal mesh point $\bz_q$ for the $j$-th coordinate is used. 
}
Accordingly, the sensing matrix $\bLambda\in\RR^{dNM\times dQ}$ is defined by 
\begin{align}
\bLambda:=
\begin{pmatrix}
\bLambda_1
\\
\vdots
\\
\bLambda_d
\end{pmatrix},
\,\text{ where }\,
\bLambda_i:=
\begin{pmatrix}
\bLambda_i^1
\\
\vdots
\\
\bLambda_i^N
\end{pmatrix}
\,\text{ with }\,
\bLambda_i^n:=
\begin{pmatrix}
\bLambda_{i,1}^{n,1} & \cdots & \bLambda_{i,d}^{n,1}
\\
\vdots & \ddots & \vdots
\\
\bLambda_{i,1}^{n,M} & \cdots & \bLambda_{i,d}^{n,M}
\end{pmatrix}.
\label{bLambda}
\end{align}
Then, the discrete version of the relationship \eqref{Gbc2}  is given by 
\begin{align*}
\Gbc\approx \bLambda\fb. 
\end{align*}

In order to facilitate the ensuing discussion, we define the wrap-around structured Hankel matrix associated to  $\fb_i\in\RR^{Q}$, for $i=1,\cdots,d$, by 
\begin{align*}
\hank_{p_i}(\fb_i):=
\begin{pmatrix}
[\fb_i]_1 & [\fb_i]_2 &\cdots & [\fb_i]_{p_i}\phantom{.1}
\\
[\fb_i]_2 & [\fb_i]_3 &\cdots & [\fb_i]_{p_i+1}
\\
\vdots & \vdots &\ddots & \vdots
\\
[\fb_i]_Q & [\fb_i]_1 &\cdots & [\fb_i]_{p_i-1}
\end{pmatrix},
\end{align*}
where $p_i<Q$ is the so-called matrix-pencil size.  
As shown in \cite{jin2016general, jin2015annihilating, ye2017deep, ye2016compressive} and reproduced in Appendix for self-containment,
 if the coordinate function $[\bF]_i$ corresponds to a smoothly varying perturbation or it has either edges or patterns,  then the corresponding Fourier spectrum $\hat{f}(\bk)$ is mostly concentrated in a small number of coefficients. 
Thus,  if  $\fb_i$ is a discretization of $[\bF]_i$ at the Nyquist sampling rate, then  according to the sampling theory of the signals with the finite rate of innovations (FRI)  \cite{vetterli2002sampling}, there exists an annihilating filter whose convolution with the image $\fb_i$ vanishes.
Furthermore, the annihilating filter size is determined by the sparsity level in the Fourier domain,
so  the  associated  Hankel structured matrix $\hank_{p_i}(\fb_i) \in \Rd^{Q\times {p_i}}$  in the image domain is low-rank if  the matrix-pencil size is chosen larger than the annihilating filter size. The interested readers are referred to Appendix or the references  \cite{jin2016general, jin2015annihilating, ye2017deep, ye2016compressive} for further details. 

 In the same way, it is expected that the block Hankel structured matrix of the  discrete source vector $\fb$, constructed as 
\begin{align*}
\hank_p(\fb)=
\begin{pmatrix}
\hank_{p_1}(\fb_1) &\cdots & \mathbf{0}
\\
\vdots & \ddots & \vdots
\\
 \mathbf{0} &\cdots &\hank_{p_d}(\fb_d) 
\end{pmatrix}
\in\RR^{dQ\times p},
\end{align*}  
is low-rank,  
{where $p=\sum_{i=1}^dp_i$. Let $r_i:=\rank(\hank_{p_i}(\fb_i))$ and $r:=\sum_{i=1}^d r_i$ where $\rank(\cdot)$ denotes the rank of a matrix. Then, a generic form of the low-rank Hankel structured constrained inverse problem can be formulated as
\begin{eqnarray}
&\ds\min_{\fb\in \Rd^{dQ}}  & \|\Gbc -\bLambda \fb\|^2
\notag
\\
&\mbox{subject to }  & \rank\left(\hank_p(\fb)\right) \leq r < p.
\label{eq:cost}
\end{eqnarray}
}

\begin{figure}[!hbt]
\centering
\includegraphics[width=0.4\linewidth]{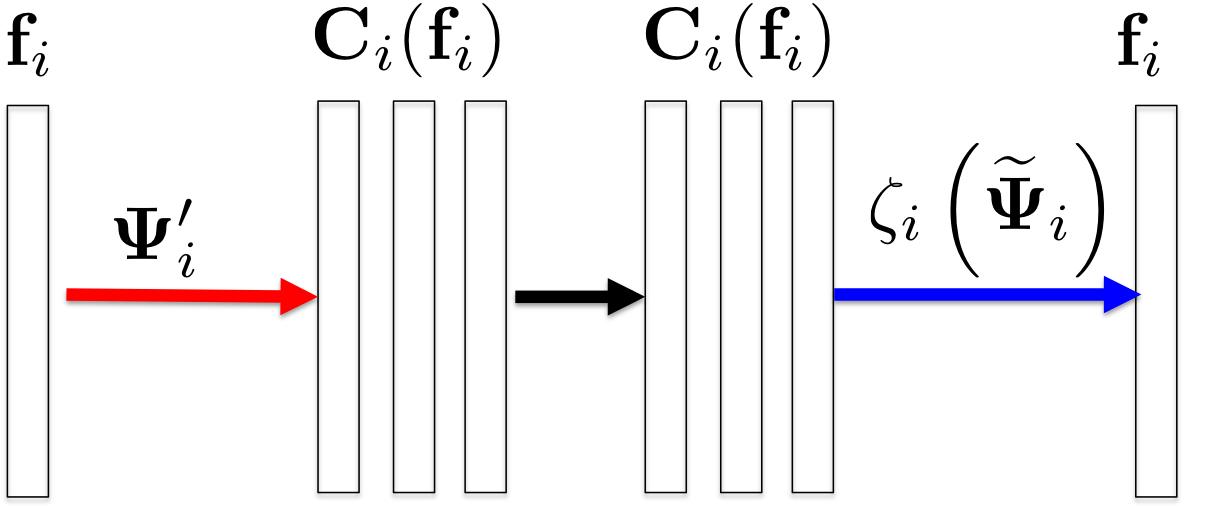}
\caption{A single layer encoder-decoder architecture in \eqref{eq:enc-dec},  
when the pooling/unpooling layers are identity matrices, i.e., $\Phib_i=\widetilde\Phib_i=\Ib_Q$.}
\label{fig:singlelayer}
\end{figure}

It is clear that, for a feasible solution $\fb=(\fb_1^\top,\cdots,\fb_d^\top)^\top$ of the regularization problem \eqref{eq:cost}, the Hankel structured matrix $\hank_{p_i}(\fb_i)$, for $i=1,\cdots,d$, admits the singular value decomposition
$\hank_{p_i}(\fb_i) = \Ub^i \Sigmab^i (\Vb^i)^{\top}$.
Here, $\Ub^i =(\bu_1^i, \cdots, \bu_{r_i}^i) \in \Rd^{Q\times r_i}$ and $\Vb^i=(\bv_1^i,\cdots, \bv_{r_i}^i)\in \Rd^{p_i\times r_i}$ denote the left and the right singular vector basis matrices, respectively, and  $\Sigmab^i=(\Sigmab^i_{kl})_{k,l=1}^{r_i}\in\mathbb{R}^{r_i\times r_i}$ refers to the diagonal matrix with singular values as elements.  
If there exist two pairs of matrices $\Phib_i, \widetilde{\Phib}_i\in\RR^{Q\times S}$ and $\Psib_i$, $\widetilde{\Psib}_i\in\RR^{p_i\times r_i}$, for each $i=1,\cdots, d$ and $S\geq Q$, satisfying the conditions
\begin{align}
\widetilde{\Phib}_i \Phib_i^\top = \Ib_{Q}
\quad\text{and}\quad 
\Psib _i\widetilde{\Psib}_i^{\top} = \bP_{\Rc(\Vb^i)}, 
\label{eq:id}
\end{align}
then 
\begin{eqnarray}
\label{eq:hankid}
\hank_{p_i}(\fb_i)  &=& \widetilde{\Phib}_i\Phib_i^\top\hank_{p_i}(\fb_i) \Psib_i \widetilde{\Psib}_i^{\top} =\widetilde{\Phib}_i \Cb_i(\fb_i)\widetilde{\Psib}_i^{\top}  = \sum_{k=1}^S\sum_{l=1}^{r_i} [\Cb_i(\fb_i)]_{kl} \widetilde{\Bb_i}^{kl} 
\end{eqnarray}
with the transformation $\Cb_i:\RR^Q\to \RR^{S\times r_i}$ given by
\begin{align}
\label{eq:coef0}
\Cb_i(\bg)= \Phib_i^\top \hank_{p_i}(\bg) \Psib_i, \qquad\forall \bg\in\RR^Q,
\end{align}
which is often called the \emph{convolutional framelet coefficient} \cite{ye2017deep}.
In Eq.~\eqref{eq:hankid}, 
\begin{eqnarray}\label{eq:hankbase}
\widetilde{\Bb_i}^{kl}:=  \widetilde\phib_{ik}\widetilde\psib_{il}^\top\in\RR^{Q\times p_i}, \quad k=1,\cdots, S,~l=1,\cdots, r_i
\end{eqnarray}
where $\widetilde\phib_{ik}$ and $\widetilde\psib_{il}$ denote the $k$-th and the $l$-th columns of $\widetilde{\Phib}_i$ and $\widetilde{\Psib}_i$, respectively.
This implies that the Hankel matrix can be decomposed using the basis matrices $\widetilde{\Bb_i}^{kl}$.
Here, the first condition in \eqref{eq:id} is the so-called \emph{frame condition}, $\Rc(\Vb^i)$ denotes the range space of $\Vb^i$ and $\bP_{\Rc(\Vb^i)}$ represents a projection onto $\Rc(\Vb^i)$ \cite{ye2017deep}.
In addition,  the pair  ($\Phib_i, \widetilde{\Phib}_i$) is \emph{non-local} in the sense that these matrices  interact with all the components of the vector $\fb_i$.  On the other hand, the pair ($\Psib_i$, $\widetilde{\Psib}_i$) is \emph{local} since these matrices interact with only $p_i$ components of $\fb_i$. 
Precisely, \eqref{eq:hankid} is equivalent to the paired encoder-decoder convolution structure when it is un-lifted to the original signal space \cite{ye2017deep}
\begin{align}
&\Cb_i(\fb_i) = \Phib_i^\top \left(\fb_i \circledast \Psib_i'\right) 
\quad\text{and}\quad 
\fb_i = 
\left(\widetilde{\Phib}_i\Cb_i(\fb_i)\right) \circledast \zeta_i\left(\widetilde{\Psib}_i\right), \label{eq:enc-dec} 
\end{align}
which is illustrated in \Cref{fig:singlelayer}.
The convolutions in \eqref{eq:enc-dec} correspond to the multi-channel convolutions (as used in standard CNN)  
with the associated filters,
\begin{align*}
\Psib_i' := 
\begin{pmatrix} 
\psib_{i1}',  \cdots,\psib_{i{r_i}}'
\end{pmatrix}\in\RR^{p_i\times r_i}
\quad\text{and}\quad
\zeta_i(\widetilde{\Psib}_i) := \frac{1}{p_i} 
\begin{pmatrix} 
\widetilde{\psib}_{i1}^\top, \cdots,\widetilde{\psib}_{i{r_i}}^\top
\end{pmatrix}^\top\in\RR^{p_ir_i}. 
\end{align*}
Here, the superposed prime over $\psib_{ik}\in\RR^{p_i}$, for fixed $i=1,\cdots,d$, and $k=1,\cdots, r_i$, indicates its flipped version, i.e., the indices of $\psib_{ik}$ are reversed \cite{ye2016compressive}.

Let us introduce the block matrices
\begin{equation*}
\begin{aligned}
&{\Phib}=\diag({\Phib}_1,\cdots, {\Phib}_d),\qquad
\widetilde{\Phib}=\diag(\widetilde{\Phib}_1,\cdots, \widetilde{\Phib}_d),
\\
&{\Psib}=\diag({\Psib}_1,\cdots, {\Psib}_d),\qquad
\widetilde{\Psib}=\diag(\widetilde{\Psib}_1,\cdots, \widetilde{\Psib}_d),
\\
&{\Vb}=\diag({\Vb}_1,\cdots, {\Vb}_d).
\end{aligned}
\end{equation*}
Then, thanks to conditions in \eqref{eq:id}, the pairs ($\Phib$, $\widetilde{\Phib}$) and ($\Psib$, $\widetilde{\Psib}$) satisfy the conditions
\begin{align*}
\widetilde{\Phib} \Phib^\top = \Ib_{dQ}
\quad\text{and}\quad 
\Psib \widetilde{\Psib}^{\top} = \bP_{\Rc(\Vb)}, 
\end{align*}
Consequently,  
\begin{align*}
\hank_p(\fb) =\widetilde{\Phib}\Phib^\top\hank_p(\fb) \Psib\widetilde{\Psib}^{\top} =\widetilde{\Phib} \Cb(\fb)\widetilde{\Psib}^{\top},
\end{align*}
with the matrix transformation $\Cb:\RR^{dQ}\to \RR^{dS\times r}$ given by
\begin{align*}
\Cb(\fb) =\diag\left(\Cb_1(\fb_1),\cdots,\Cb_d(\fb_d)\right)= \Phib^\top \hank_p(\fb) \Psib.
\end{align*}
Let $\Hbc_{i}$, for  $i=1,\cdots,d$, refer  to the space of signals admissible in the form \eqref{eq:enc-dec}, i.e., 
\begin{align*}
\Hbc_{i}:=\Big\{\bg_i 
\in\RR^{Q}\,\Big|\,\,  
\bg_i=
\left(\widetilde{\Phib}_i\Cb_i\right) \circledast \zeta_i\left(\widetilde{\Psib}_i\right),
\Cb_i(\bg_i ) = \Phib^\top_i \left(\bg_i \circledast \Psib_i'\right),\Big\}.
\end{align*}
Then,  the problem \eqref{eq:cost} can be converted to
\begin{eqnarray}
\label{eq:newcostBlock}
\min_{\fb\in\prod_{i=1}^d \Hbc_{i}} \,\left\|\Gbc -\bLambda \fb\right\|^2,
\end{eqnarray}
or equivalently,
\begin{eqnarray}
\label{eq:newcost}
\min_{\fb_i\in\Hbc_{i}} \, \left\|\Gbc_i -\bLambda_i\fb\right\|^2,\quad i=1,\cdots, d,
\end{eqnarray}
where the sub-matrices $\Gbc_i$ and  $\bLambda_i$ are defined in  \cref{Gbc,bLambda}, respectively. 

\subsection{Extension to Deep Neural Network}\label{sect:CNN}


One of the most important discoveries in \cite{ye2017deep} is that an encoder-decoder network architecture in the convolutional neural network (CNN) is emerged from Eqs.~\eqref{eq:hankid},\eqref{eq:coef0}, and  \eqref{eq:enc-dec}. 
In particular, the non-local bases matrices  $\Phib_i$ and $\widetilde{\Phib}_i$ play the role of user-specified \emph{pooling} and \emph{unpooling} operations, respectively  (see,  \cref{ss:DUN}),
whereas the local-bases $\Psib_i$ and  $\widetilde{\Psib}_i$ correspond  to the \emph{ encoder and decoder layer convolutional filters} that have to be learned from the data
\cite{ye2017deep}.

Specifically, our goal is to learn ($\Psib_i,\widetilde{\Psib}_i$) in a \emph{data-driven} fashion so that the optimization problem \eqref{eq:newcostBlock} (or equivalently \eqref{eq:newcost}) can be simplified. Toward this, we first define $\bLambda^\dagger$ (resp. $\bLambda^\dagger_i$) as a right pseudo-inverse of $\bLambda$ (resp. $\bLambda_i$), i.e., $\bLambda\bLambda^\dagger \Gbc =\Gbc$ (resp. $\bLambda_i\bLambda_i^\dagger \Gbc_i =\Gbc_i$) for all $\Gbc\in\RR^{dMN}$ so that the cost in \eqref{eq:newcostBlock} (resp \eqref{eq:newcost}) can be automatically minimized with
the right pseudo-inverse solution. However,  the solution leads to 
\begin{align*}
\fb=\bLambda^\dagger \Gbc=  \fb^*+\fb^0:=\begin{pmatrix} \fb^*_1\\ \vdots \\\fb^*_d \end{pmatrix}+\begin{pmatrix} \fb^0_1\\ \vdots \\\fb^0_d \end{pmatrix},
\end{align*}
where $\fb^*$ denotes the true solution and  $\fb^0\in \Nc(\bLambda)$.  
{
Therefore, one looks for the matrices ($\Psib_i,\widetilde{\Psib}_i$)  such that
\begin{align*}
\fb^*_i =&  
\Kbc_i[\Psib_i,\widetilde{\Psib}_i](\fb_i),
\qquad\forall\,\, i=1,\cdots, d,
\end{align*}
where the operator $\Kbc_i[\Psib_i,\widetilde{\Psib}_i]:\RR^Q\to \RR^Q$ is defined in terms of the mapping ${\Cb}_i(\cdot)$ = ${\Cb}_i[\Psib_i](\cdot)$ as 
\begin{align}
\Kbc_i[\Psib_i,\widetilde{\Psib}_i](\fb_i)=\Kbc_i[\Psib_i,\widetilde{\Psib}_i](\fb_i^*+\fb_i^0)
:= \left(\widetilde{\Phib}_i {\Cb}_i[\Psib_i](\fb^*_i+\fb_i^0)\right) \circledast \zeta_i(\widetilde{\Psib}_i), \label{eq:Ki}
\end{align}
for all $\fb=\fb^*\oplus\fb^0\in \Rc(\bLambda^\dagger)\oplus\Nc(\bLambda)$. In fact, the operator $\Kbc_i[\Psib_i,\widetilde{\Psib}_i]$ in \eqref{eq:Ki} can be engineered so that its output $\fb$ belongs to $\Rc(\bLambda^\dagger)$. This can be achieved by selecting the filters $\Psib'_i$'s that \emph{annihilate} the null-space components $\fb_i^0$'s, i.e., 
\begin{align*}
\fb_i^0 \circledast \Psib'_i\approx \mathbf{0},\quad i=1,\cdots, d,
\end{align*} 
so that 
\begin{align*}
{\Cb}_i[\Psib_i](\fb^*_i+\fb_i^0)=\Phib_i^\top\left(\fb_i^*\oplus\fb_i^0 \circledast \Psib'_i\right)\approx\Phib_i^\top\left(\fb_i^*\circledast \Psib'_i\right),\qquad\forall i=1,\cdots, d,
\end{align*}
In other words, the block filter $\Psib'$ should span the orthogonal complement of  $\Nc(\bLambda)$.}
Therefore,  the local bases learning problem becomes
\begin{align}
\label{eq:opt}
\min_{(\Psib_i,\widetilde{\Psib}_i)} \sum_{\ell=1}^L\left\|\fb_i^{*(\ell)} - \Kbc_i[\Psib_i,\widetilde{\Psib}_i](\fb_i^{(\ell)})\right\|^2, \quad i=1,\cdots, d,
\end{align}
where 
$$
\left\{\left(\fb^{(\ell)}_i,\fb^{*(\ell)}_i\right):=\left(\left[\bLambda^\dagger \Gbc^{(\ell)}\right]_i,\fb^{*(\ell)}_i\right)\right\}_{\ell=1}^L,
$$ 
denotes the training data-set composed of the input and ground-truth pairs. 
This is equivalent to saying that the proposed neural network is for learning the local basis from the training data assuming that the Hankel matrices associated with discrete source densities $\fb_i$ are of rank $r_i$ \cite{ye2017deep}.

Still, the convolutional framelet expansion is  linear, so  we restricted the space so that
the framelet coefficient matrices $\Cb_i(\fb_i)$ are restricted to have positive elements only, i.e.,
 the signal lives in the conic hull of the  basis to enable part-by-part representation similar to non-negative matrix factorization (NMF) \cite{
lee1997unsupervised, lee1999learning,  lee2001algorithms}: 
 \begin{align*}
\Hbc_{i}^0:=\Big\{\bg
\in\RR^{Q}\,\Big|\,\,  \bg &=
\left(\widetilde{\Phib}_i\Cb_i(\bg)\right) \circledast \zeta_i\left(\widetilde{\Psib}_i\right),
\nonumber
\\
&\Cb_i (\bg)= \Phib^\top_i \left(\bg \circledast \Psib_i'\right), \quad [\Cb_i(\bg)]_{kl}\geq 0, \,\,\forall k,l\Big\},
\end{align*}
for $i=1,\cdots, d$. 
This positivity constraint can be implemented using the \emph{rectified linear unit} (ReLU) \cite{nair2010rectified} during training.
Accordingly, the local basis learning problem \eqref{eq:opt} can equivalently be expressed as
\begin{align*}
\min_{(\Psib_i,\widetilde{\Psib}_i)} \sum_{\ell=1}^L\left\|\fb_i^{*(\ell)} - \Kc^\varrho_i[\Psib_i,\widetilde{\Psib}_i]\left(\fb_i^{(\ell)}\right)\right\|^2, \quad i=1,\cdots, d.
\end{align*}
Here, the operator $\Kbc_i^\varrho[\Psib_i,\widetilde{\Psib}_i]:\RR^Q\to\RR^Q$ is defined analogously as in \eqref{eq:Ki} but in terms of the mapping $\Cb_i^\varrho:\RR^Q\to\RR^{S\times r_i}$ given by  
\begin{align*}
\Cb_i^\varrho(\bg) =\varrho\left(\Phib_i^\top \left(\bg \circledast \Psib_i'\right)\right), \quad \forall\,\, i=1,\cdots, d,
\end{align*}
where $\varrho$ denotes ReLU, i.e., for arbitrary matrix $\mathbf{A}\in\RR^{S\times r_i}$, we have $[\varrho(\mathbf{A})]_{kl}\geq 0$, for all $k$ and $l$.

\begin{figure}[!hbt]
\centering
\includegraphics[width=0.7\linewidth]{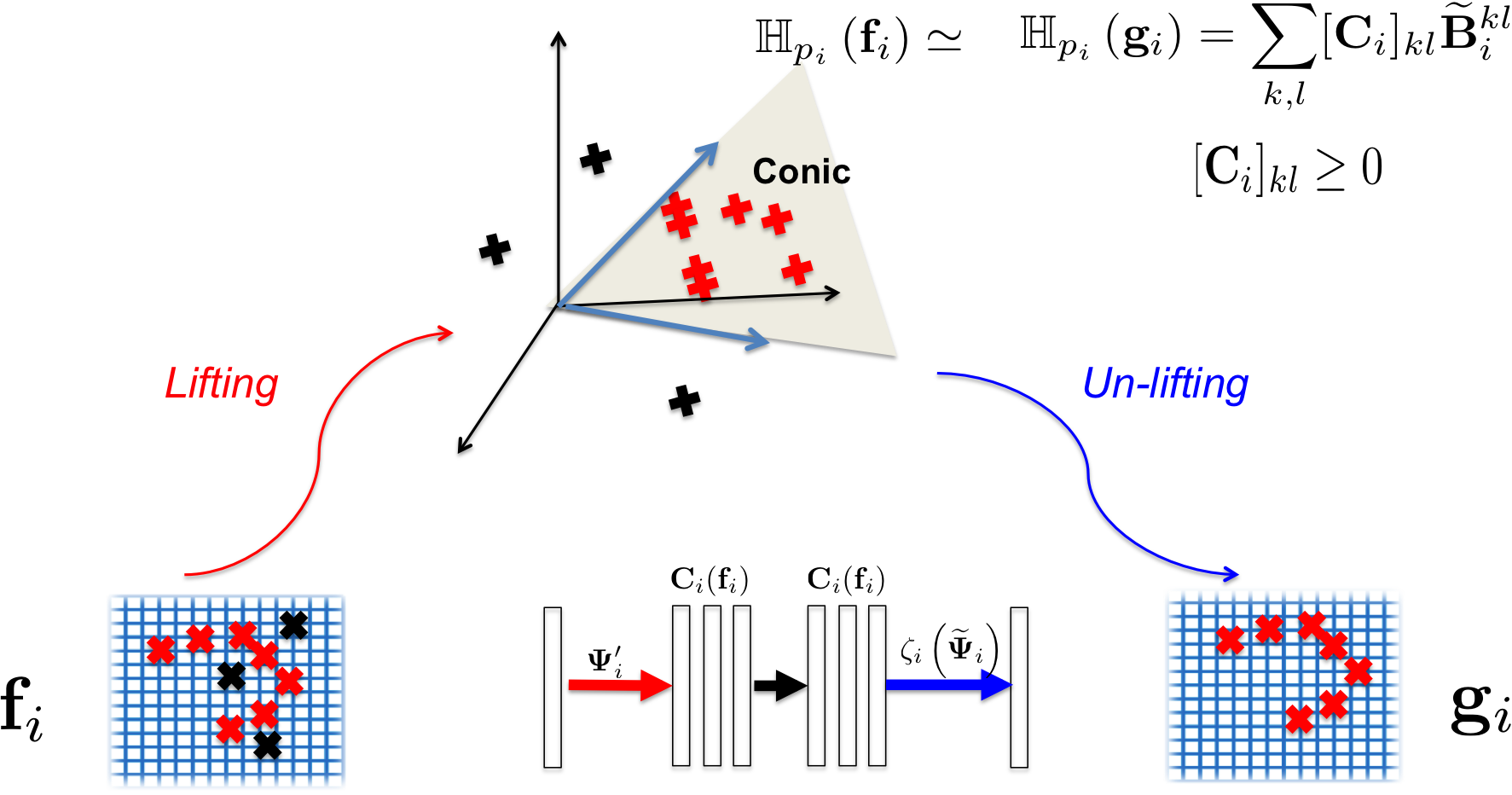}
\caption{Geometry of single layer encoder decoder network for denoising. The signal is first lifted into higher dimensional space, which is then decomposed into the positive combination of
bases. During this procedure, the outlier signals (black color) are placed outside of the conic hull of the bases, so that they can be removed during the decomposition. When this high conic decomposition procedure is observed in the original signal space, it becomes one level encoder-decoder neural network with ReLU.
}
\label{fig:geometry}
\end{figure}

The geometric implication of this representation is illustrated in
\Cref{fig:geometry}.
Specifically,  the original image $\fb_i$ is  first \emph{lifted} to higher dimensional space via Hankel matrix, $\hank_{p_i}(\fb_i)$,
which is then decomposed into positive (conic) combination using the matrix bases $\widetilde\Bb_i^{kl}$ in \eqref{eq:hankbase}.
During this procedure, the outlier signals (black color) are placed outside of the conic hull of the bases, so that they can be removed during the decomposition. When this high conic decomposition procedure
is observed in the original signal space, it becomes one level encoder-decoder neural network with ReLU.  Therefore, an encoder-decoder network can be understood as a signal space manifestation of the conic decomposition of the signal being lifted to  a higher-dimensional space.

\begin{figure}[!hbt]
\centering
\includegraphics[width=0.7\linewidth]{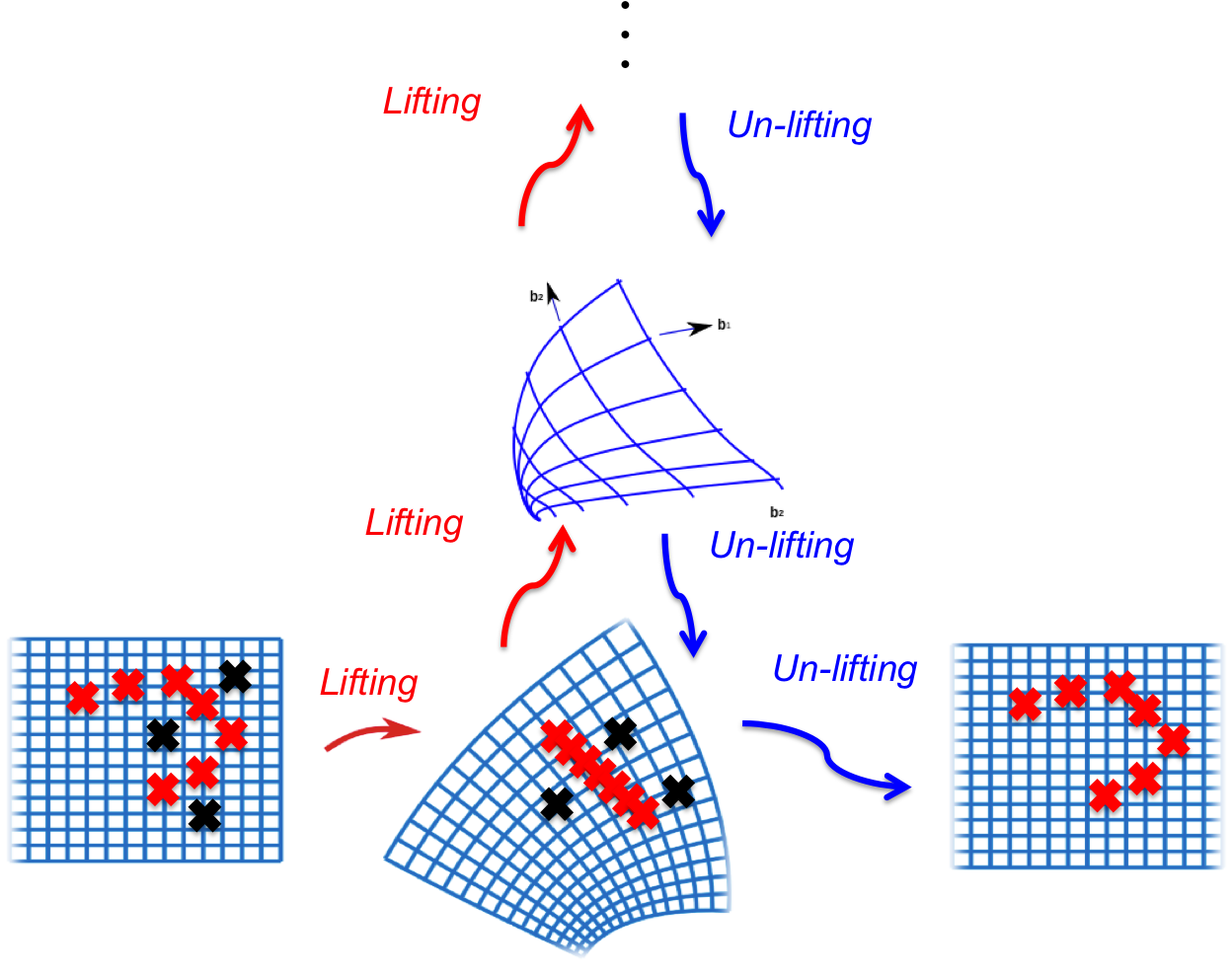}
\centerline{\mbox{(a)}}
\includegraphics[width=0.7\linewidth]{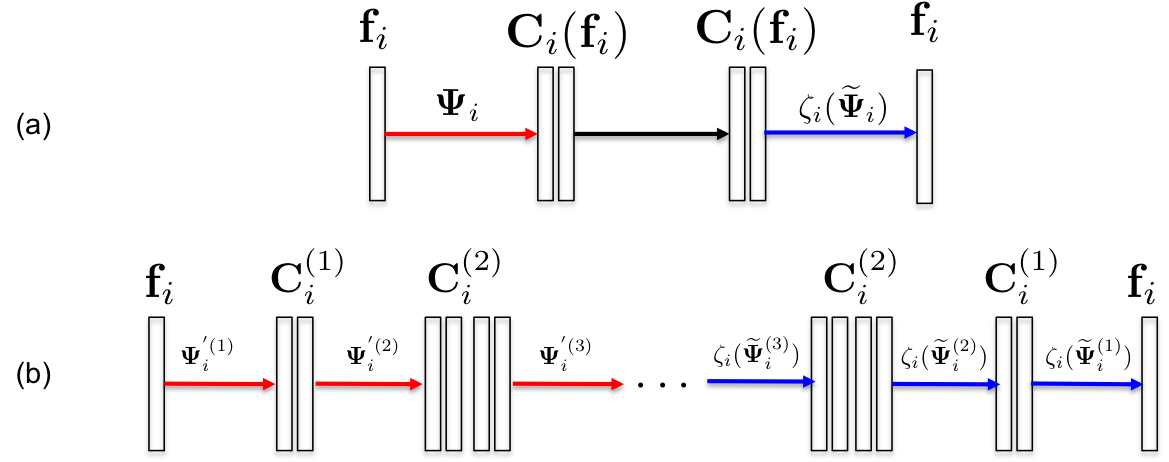}
\centerline{\mbox{(b)}}
\caption{(a)Geometry of multi-layer encoder decoder network, and (b) its original space manifestation as a multi-layer encoder-decoder network.
}
\label{fig:multilayer}
\end{figure}

The idea can be further extended to the multi-layer deep neural network.
Specially, suppose that the encoder and decoder convolution filter $\Psib_i'$ and $\zeta_i(\widetilde\Psib_i)$
can be represented in a cascaded convolution of small length filters:
\begin{eqnarray*}
\Psib'_i &=&\Psib_i^{'(0)} \circledast \cdots \circledast \overline \Psib^{'(J)} \\
\zeta_i(\widetilde\Psib_i) &=& \zeta_i(\widetilde \Psib^{(J)}) \circledast \cdots \circledast  \zeta_i(\widetilde\Psib^{(0)}), 
\end{eqnarray*}
then the signal space is recursively defined as 
 \begin{align*}
\Hbc_{i}^0:=\Big\{\bg
\in\RR^{Q}\,\Big|\,\, &  \bg=
\left(\widetilde{\Phib}_i\Cb_i(\bg)\right) \circledast \zeta_i\left(\widetilde{\Psib}_i\right),
\nonumber
\\
&\Cb_i (\bg)= \Phib^\top_i \left(\bg \circledast \Psib_i'\right)\in\Hbc_{i}^1, \quad [\Cb_i(\bg)]_{kl}\geq 0, \,\,\forall k,l\Big\},
\end{align*}
where, for all $\jmath=1,\cdots, J-1\in\NN$,
\begin{align}
&\Hbc_{i}^\jmath:=\Big\{\bA
\in\RR^{Q\times Q_i^{(\jmath)}}\,\Big|\,\, \bA=
\left(\widetilde{\Phib}_i\Cb_i^{(\jmath)}(\bA)\right) \circledast \zeta_i\left(\widetilde{\Psib}_i^{(\jmath)}\right),
\nonumber
\\
&\qquad\qquad\qquad\Cb_i^{(\jmath)} (\bA)= \Phib^\top_i \left(\bA \circledast \Psib_i'^{(\jmath)}\right)\in\Hbc_{i}^{\jmath+1}, \,\, [\Cb_i(\bg)]_{kl}\geq 0, \,\,\forall k,l \Big\},  \label{eq:multilayer}
\\
&\Hbc_{i}^J:= \RR^{Q\times Q_i^{(J)}}. \notag
\end{align}
Here, the $\jmath$-th layer encoder  and decoder filters, $\Psib_i'^{(\jmath)}\in\RR^{p_i^{(\jmath)}Q_i^{(\jmath)}\times R_i^{(\jmath)}}$ and $\zeta_i\left(\widetilde{\Psib}_i^{(\jmath)}\right)\in\RR^{p_i^{(\jmath)}R_i^{(\jmath)}\times Q_i^{(\jmath)}}$, are given  by 
\begin{align*}
&\Psib_i'^{(\jmath)}:=
\begin{pmatrix}
\psi'^{1}_{i1} & \cdots & \psi'^{1}_{iR_i^{(\jmath)}}
\\
\vdots & \ddots & \vdots
\\
\psi'^{Q^{(\jmath)}_i}_{i1} & \cdots & \psi'^{Q^{(\jmath)}_i}_{iR_i^{(\jmath)}}
\end{pmatrix}
\,\text{ and }\,
\zeta_i\left(\widetilde{\Psib}_i^{(\jmath)}\right):=
\begin{pmatrix}
\widetilde{\psi}^{1}_{i1} & \cdots & \widetilde{\psi}^{Q_i^{(\jmath)}}_{i1}
\\
\vdots & \ddots & \vdots
\\
\widetilde{\psi}^{1}_{iR^{(\jmath)}_i} & \cdots & \widetilde{\psi}^{Q_i^{(\jmath)}}_{iR^{(\jmath)}_i}
\end{pmatrix},
\end{align*}
where $p_i^{(\jmath)}$, $Q_i^{(\jmath)}$, and $R_i^{(\jmath)}$ are the filter lengths, the number of input channels, and the number of output channels, respectively.
This is equivalent to recursively applying high dimensional conic decomposition procedure  to the next level
convolutional framelet coefficients as illustrated in \Cref{fig:multilayer}(a).  The resulting signal space manifestation 
is a deep neural network shown  in \Cref{fig:multilayer}(b).

\subsection{Dual-Frame U-Net}\label{ss:DUN}

\begin{figure}[!hbt]
\centering
\includegraphics[width=0.7\linewidth]{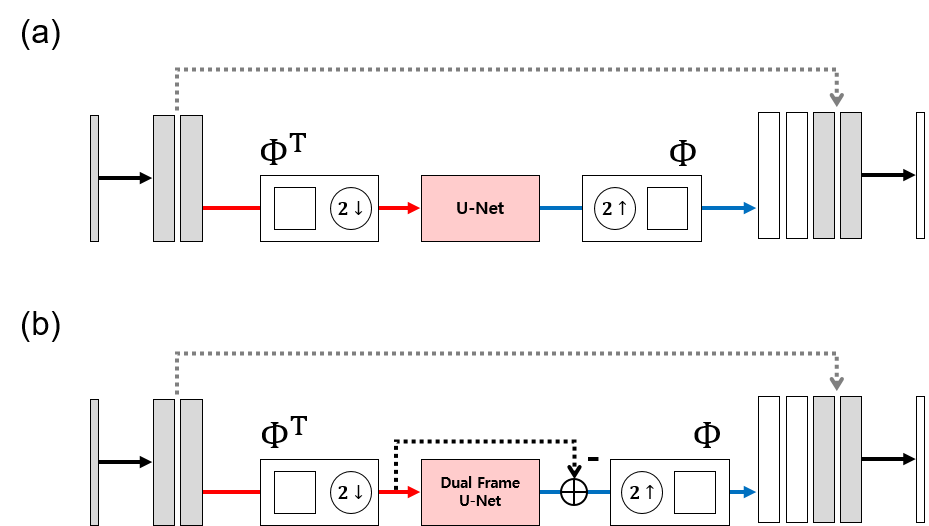}
\caption{Simplified U-Net architecture and its variant. (a) Standard U-Net and
(b) Dual-Frame U-Net.}
\label{fig:netArch}
\end{figure}

As discussed before, 
the non-local bases $ \Phib_i^\top$ and $\widetilde{\Phib}_i$ correspond  to the generalized pooling and unpooling operations, which can be designed by the users for specific inverse problems. Here,  the key requisite is the frame condition  in \eqref{eq:id}, i.e., $\widetilde{\Phib}_i \Phib_i^\top = \Ib_{Q}$. As the artifacts of the time-reversal recovery from sparse measurements are distributed globally, a network architecture with large receptive fields is needed. Thus, in order to learn the optimal local basis from the minimization problem \eqref{eq:opt}, we adopt the commonly used CNN architecture known as \emph{U-Net} \cite{ronneberger2015u} and its deep convolutional framelets based variant, coined as \emph{Dual-Frame U-Net} \cite{han2017framing}  (see \Cref{fig:netArch}).  These networks have pooling layers with down sampling, resulting exponentially large receptive fields.

As shown in  \cite{han2017framing},  one of the main limitation of the standard U-Net is that it does not satisfy the frame condition in \eqref{eq:id}.  Specifically, by considering both skipped connection and the pooling $\Phib_i^\top$ in \Cref{fig:netArch}(a),
the non-local basis for the standard U-Net is given by
\begin{eqnarray}
\label{eq:Pext}
\Phib_{i}^\top := \begin{pmatrix}  {\Ib}_{Q}  \\ \Phib_{i,{\rm avg}}^\top \end{pmatrix}\in\RR^{\frac{3Q}{2}\times Q}, 
\end{eqnarray}
where 
$\Phib_{i,{\rm avg}}^\top$ denotes an average pooling operator given by
\begin{eqnarray*}
\Phib_{i,{\rm avg}}^\top  =  
 \frac{1}{\sqrt{2}}  
 \begin{pmatrix} 
 1 & 1 &  0 & 0 & \cdots & 0  & 0  \\  
 0 & 0 & 1 & 1 & \cdots & 0  & 0  \\
 \vdots & \vdots & \vdots  & \vdots & \ddots &  & \vdots  \\
  0 & 0 &  0 & 0 & \cdots & 1 & 1 
  \end{pmatrix} 
  \in \RR^{\frac{Q}{2}\times Q}.
\end{eqnarray*}
Moreover, the unpooling layer in the standard U-Net is given by  ${\Phib}_i$. Therefore, 
\begin{align*}
\widetilde{\Phib}_i \Phib_i^\top =  {\Ib}_{Q} +\Phib_{i,{\rm avg}}\Phib_{i,{\rm avg}}^\top \neq \Ib_Q.
\end{align*}
Consequently, the frame condition in \eqref{eq:id} is not satisfied.
As shown in \cite{ye2017deep}, this results in the duplication of the low-frequency components, making the final results blurry.

In order to address the aforementioned problem, in the Dual-Frame U-Net \cite{han2017framing}, the dual 
frame is directly implemented. More specifically, the dual frame $\widetilde{\Phib}_i$ for the specific frame operator \eqref{eq:Pext}
is given by
\begin{eqnarray*}
\widetilde{\Phib}_{i} :=  (\Phib_{i}\Phib_{i}^\top)^{-1}\Phib_{i} &=&  ( \Ib_Q + \Phib_{i,{\rm avg}} \Phib_{i,{\rm avg}}^\top)^{-1} 
\begin{bmatrix}  {\Ib}_Q  & \Phib_{i,{\rm avg}} \end{bmatrix}  
\notag
\\
&=& 
\begin{bmatrix} 
\Ib_Q-\Phib_{i,{\rm avg}}\Phib_{i,{\rm avg}}^\top/2 & \Phib_{i,{\rm avg}}/2 
\end{bmatrix},
\end{eqnarray*}
where the matrix inversion lemma and the orthogonality $\Phib_{i,{\rm avg}}^\top\Phib_{i,{\rm avg}} =\Ib_Q$ are invoked to arrive at the last equality.
It was shown in  \cite{han2017framing} that the corresponding generalized unpooling operation is given by 
\begin{eqnarray} 
\label{eq:dualUnet}
\widetilde{\Phib}_i \begin{bmatrix} \Bb_i  \\ \Cb_i(\gb)   \end{bmatrix}
&=& \Bb_i - \frac{1}{2} \underbrace{\Phib_{i,{\rm avg}}}_{\mbox{unpooling}} \overbrace{(\Phib_{i,{\rm avg}}^\top \Bb_i-  \Cb_i(\gb))}^{\mbox{residual}}, 
\end{eqnarray}
where $\Bb_i$ denotes the skipped component.
Equation \eqref{eq:dualUnet} suggests a network structure for the Dual-Frame U-Net.
More specifically,  unlike the U-Net,  the \emph{residual signal} at the low resolution  should be upsampled through the unpooling layer and subtracted from the by-pass signal to eliminate the duplicate contribution of the low-frequency components.
This can be easily implemented using additional bypass connection for the low-resolution signal as shown in \Cref{fig:netArch}(b).  
This simple fix allows the proposed network to satisfy the frame condition \eqref{eq:id}. 
The interested readers are suggested to consult \cite{han2017framing} for further details.

\section{Network design and training}\label{Sect:NT}

Let us now design the U-Net and Dual-Frame U-Net neural networks for the elastic source imaging based on the analysis performed in the previous section. For simplicity,  consider a 2D  case (i.e., $d=2$) for the recovery of the $x$-component (i.e., $[\bF]_1$) of the unknown source. The $y$-component (i.e., $[\bF]_2$) of the source can be obtained in exactly the same fashion using the  network architectures  discussed above. 

\subsection{Description of the forward solver and time-reversal algorithm}

For numerical illustrations and generation of the training data, the region of interest $\Omega$ is considered as a unit disk centered at origin. Each solution of the elastic wave equation is computed over the box $B = \big[-\beta/2,\beta/2\big]^2$ so that $\Omega \subset B$, i.e.,  $(\bx,t) \in\big[-\beta/2, \beta/2\big]^2 \times\big[0,\T\big]$ with $\beta = 4$ and $\T=2$.  The temporal and spatial discretization steps are, respectively, chosen to be  $h_t = 2^{-6}\T $ and $h_x =2^{-7}\beta $. The Lam\'e parameters are chosen in such a way that the pressure and the shear wave speeds in the medium are, respectively,  $c_P=\sqrt{3}m.s^{-1}$ and $c_S=1m.s^{-1}$.  

The  Lam\'e system
\begin{equation*}
\begin{cases}
\ds\frac{\partial^2 \bu}{\partial t^2 }(\bx,t) - \LE\bu(\bx,t)=\frac{\partial \delta_0}{\partial t}\bF(\bx),
\quad \quad (\bx ,t) \in \RR^2 \times \RR , 
\\
\ds\bu(\bx,0) =\mathbf{0}  \quad \text{and} \quad   \ds 
\frac{\partial \bu}{\partial t}(\bx,0) = {\bf 0},
\end{cases}
\end{equation*}
is numerically solved over the  box $B$ with periodic boundary conditions. A \emph{splitting spectral Fourier}  approach \cite{spectral} is used together with a perfectly matched layer (PML) technique \cite{pml} to simulate a free outgoing interface on $\partial B$.  The weighted time-reversal function  $\I_{\rm WTR}(\bx)$ also requires a Helmholtz decomposition algorithm. Since the support of the function $\I_{\rm WTR}(\bx)$  is included in $\Omega \subset B$, a Neumann boundary condition is used on $\partial B$ and a weak Neumann problem is solved in order to derive the Helmholtz decomposition. This decomposition is numerically obtained with a fast algorithm proposed in \cite{Wiegmann} based on a symmetry principle and a Fourier Helmholtz decomposition algorithm.  The interested readers are suggested to consult \cite[Sect. 2.2.1]{TrElastic} for more details on the numerical algorithm.  

\subsection{Data preparation}\label{sect:DataPrep}

As a training data-set,  training pairs $\{(\fb_i^{(\ell)}, \fb_i^{*(\ell)})\}_{\ell=1}^L$ are generated with $L = 5000$ where $\fb_i^{(\ell)}$ is a numerically generated input image  and  $\fb_i^{*(\ell)}$ denotes  the synthetic ground-truth phantoms. More specifically, the input images  $\{\fb_i^{(\ell)}\}_{\ell=1}^L$ are generated numerically by first computing the solution formula for the wave equation for a set of phantom images $\{\fb_i^{*(\ell)}\}_{\ell=1}^L$  and then applying the time-reversal algorithm. The pixel values of input images are centered at origin by subtracting the mean intensity of each individual image and dividing it by the maximum value over the entire data-set. The phantoms are generated using the in-built MATLAB  \emph{phantom} function such that each phantom had up to ten random overlapping ellipses with their supports compactly contained in $\Omega$. The centers of the ellipses are randomly selected from $[-0.375,0.375]$. The minor and major axes are chosen as random numbers from $[-0.525, 0.525]$. The angles between the horizontal semi-axes of the ellipses and the $x$-axis of the image are also randomly selected from $[-\pi,\pi]$. The intensity values of the ellipses are restricted between $[-10,10]$ so that the values of the overlapping area are negatively or positively added. Finally,  every generated phantom is normalized by subtracting the minimum value and dividing the maximum value sequentially so that its intensity lies in a positive range $[0,1]$.

\subsection{Network architectures}

\begin{figure}[!hbt]
\centering
\includegraphics[width=0.9\linewidth]{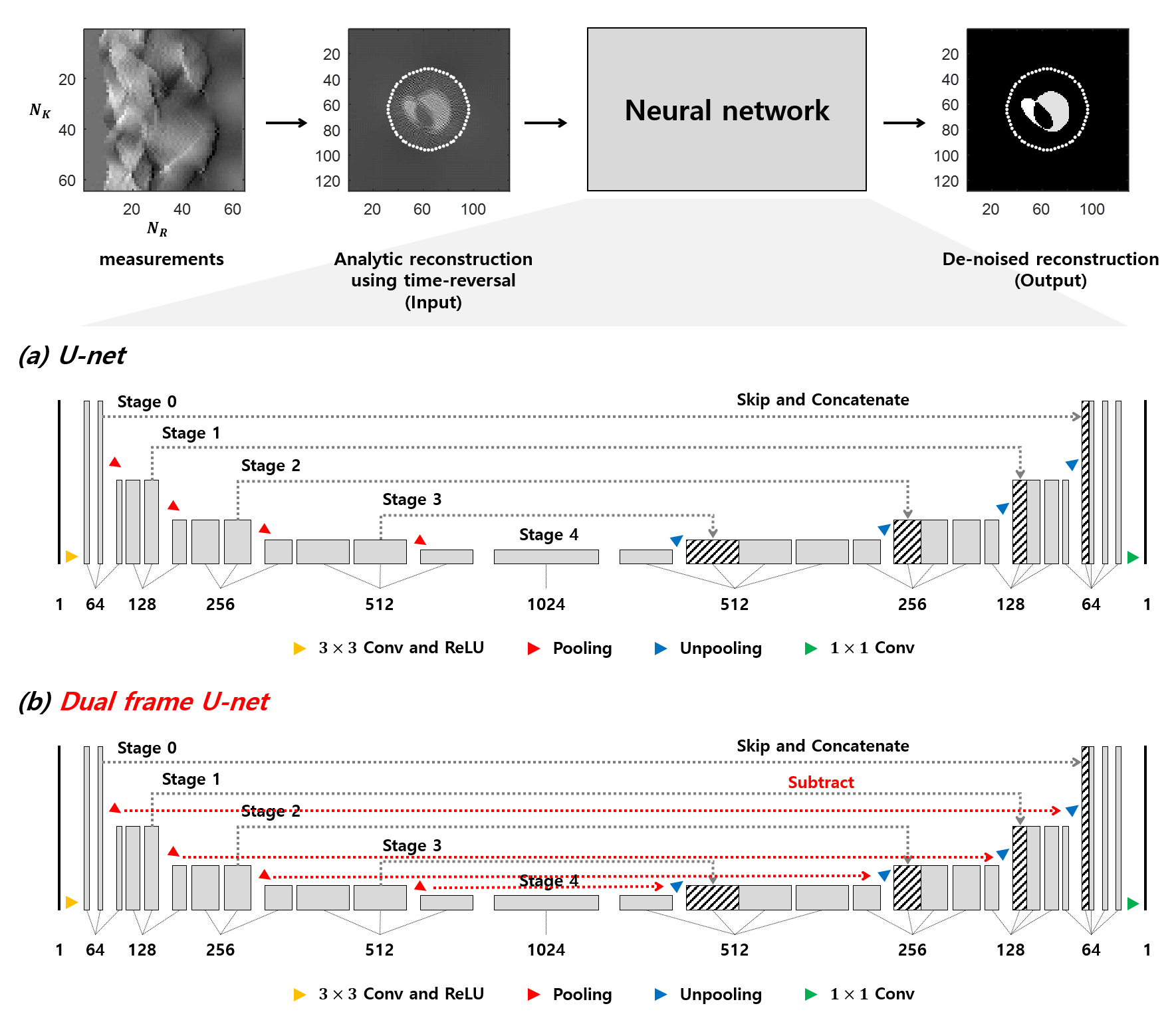}
\caption{Schematic illustration}
\label{fig:scheme}
\end{figure}

The original and Dual-Frame U-Nets consist of convolution layer, ReLU, and contracting path connection with concatenation (\Cref{fig:scheme}). Specifically, each stage contains four sequential layers composed of convolution with $3\times 3$ kernels and ReLU layers. Finally, the last stage has two sequential layers and the last layer contains only a single convolution layer with $1\times 1$ kernel. The number of channels for each convolution layer is illustrated in \Cref{fig:scheme}. Note that the number of channels is doubled after each max pooling layer. The differences between the original and Dual-Frame U-Nets are from additional residual paths illustrated in \Cref{fig:netArch}.

\subsection{Network training}

In order to deal with the imbalanced distribution of non-zero values in the label phantom images $\fb_i^{*(\ell)}$ and to prevent the proposed network to learning a trivial mapping (rendering all zero values), the non-zero values are weighted by multiplying a constant according to the ratio of the total number of voxel over the non-zero voxels. All the convolutional layers were preceded by appropriate zero-padding to preserve the size of the input. The mean squared error (MSE) is used as a loss function and the network is implemented using {Keras library} \cite{chollet2015keras}. The weights for all the convolutional layers were initialized using {Xavier initialization}. The generated data  is divided into $4000$ training and $1000$ validation data-sets. 
For training,  the batch size of $64$ and  {Adam optimizer} \cite{kingma2014adam} with the default parameters as mentioned in the original paper are used, i.e., the learning rate $=0.0001$, $\beta_1 =0.9$, and $\beta_2 =0.999$ are adopted. The training runs for up to $200$ epochs with early stopping if the validation loss has not improved in the last $20$ epochs.  GTX 1080 graphic processor and i7-6700 CPU ($3.40$ GHz) are used. The network took approximately $1300$ seconds. 

\section{Numerical experiments and discussion}\label{Sect:NE}
 
In this section,  some numerical realizations of the proposed algorithm are presented for the resolution of the inverse elastic source problem and the performances of the proposed deep learning frameworks are debated.  The examples of sparse targets with binary intensities  and extended targets with variable intensities are discussed. The sparse targets are modeled by an elongated tubular shape and an ellipse. The extended targets are modeled by the Shepp-Logan phantom. The performance of the proposed framework is compared with the results rendered by the  weighted time-reversal algorithm with sub-sampled sparse data and total variation (TV) based regularization approach applied on the low-resolution images provided by the time-reversal algorithm. The reconstructed images are compared under both clean and noisy measurement conditions with the TV- regulatization using {fast iterative shrinkage threshholding algorithm} (FISTA) of Beck and Teboulle \cite{beck2009fast}.  For comparison,   the {peak-signal-to-noise ratio} (PSNR) and the {structural similarity index} (SSIM) are  used as metrics, where
\begin{align*}
& {\rm PSNR}:=20\,\log_{10}\left(\frac{\widetilde{N}\widetilde{M}\|\hat{\fb}_1\|^2_\infty}{\|\hat{\fb}_1-\fb^*_1\|_2}\right),
\quad
{\rm SSIM}:=\frac{\left(2\mu_{\hat{\fb}_1}\mu_{\fb^*_1}+c_1\right)\left(2\sigma_{\hat{\fb}_1\fb^*_1}+c_2\right)}{\left(\mu_{\hat{\fb}_1}^2+\mu_{\fb^*_1}^2+c_1\right)\left(\sigma_{\hat{\fb}_1}^2+\sigma^2_{\fb^*_1}+c_2\right)}.
\end{align*}
Here, $\widetilde{M}$ and $\widetilde{N}$ are the number of pixels in the rows and columns, $\hat{\fb}_1$ and $\fb^*_1$ are the reconstructed image and ground truth, $\mu_{\hat{\fb}_1}$ and $\mu_{\fb^*_1}$ are the expectations, $\sigma^2_{\hat{\fb}_1}$ and $\sigma^2_{\fb^*_1}$ are the variances,  and $\sigma^2_{\hat{\fb}_1\fb^*_1}$ is the covariance of $\hat{\fb}_1$ and $\fb^*_1$, respectively. Here, $c_1$ and $c_2$ are stabilization parameters and are chosen as $c_1=(0.01\xi)^2$ and $c_1=(0.03\xi)^2$ with $\xi$ being the dynamic range of the pixel intensity.

\subsection{Results}
\begin{figure}[!hbt]
\centering
\includegraphics[width=0.8\linewidth]{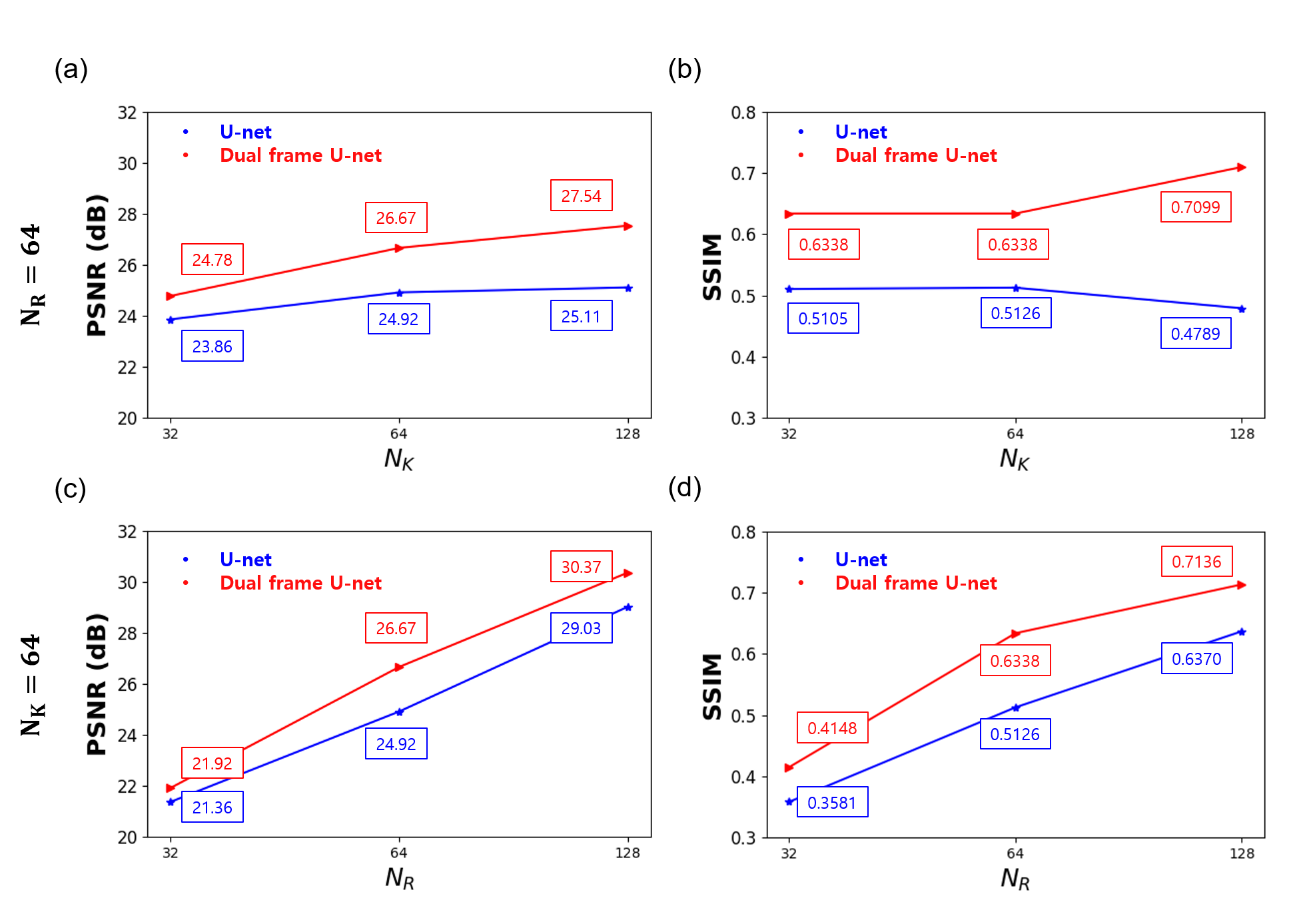}
\caption{PSNR (left column) and SSIM (right column) results of the standard U-Net (blue) and Dual-Frame U-Net (red).}
\label{fig:psnr_ssim}
\end{figure}

\Cref{fig:psnr_ssim} shows the variations of PSNR and SSIM values of the reconstructed test images using the standard and Dual-Frame U-Net. By increasing the number of recorders ($N_R$) or scanning rate ($N_K$), it is observed that the PSNR and SSIM values of the images show a monotonically increasing trend except  for the SSIM value of the image from U-Net with $N_K=128$ (\Cref{fig:psnr_ssim}(b)). 
On the other hand, the performance of the Dual-Frame U-Net always improved with more measurement data. In addition, the PSNR and SSIM values of the images from the Dual-Frame U-Net are always higher than the ones from the standard U-Net. This  suggests that the Dual-Frame U-Net, which satisfies the frame condition, is a robust and predictable reconstruction scheme.

\begin{figure}[!hbt]
\centering
\includegraphics[width=0.8\linewidth]{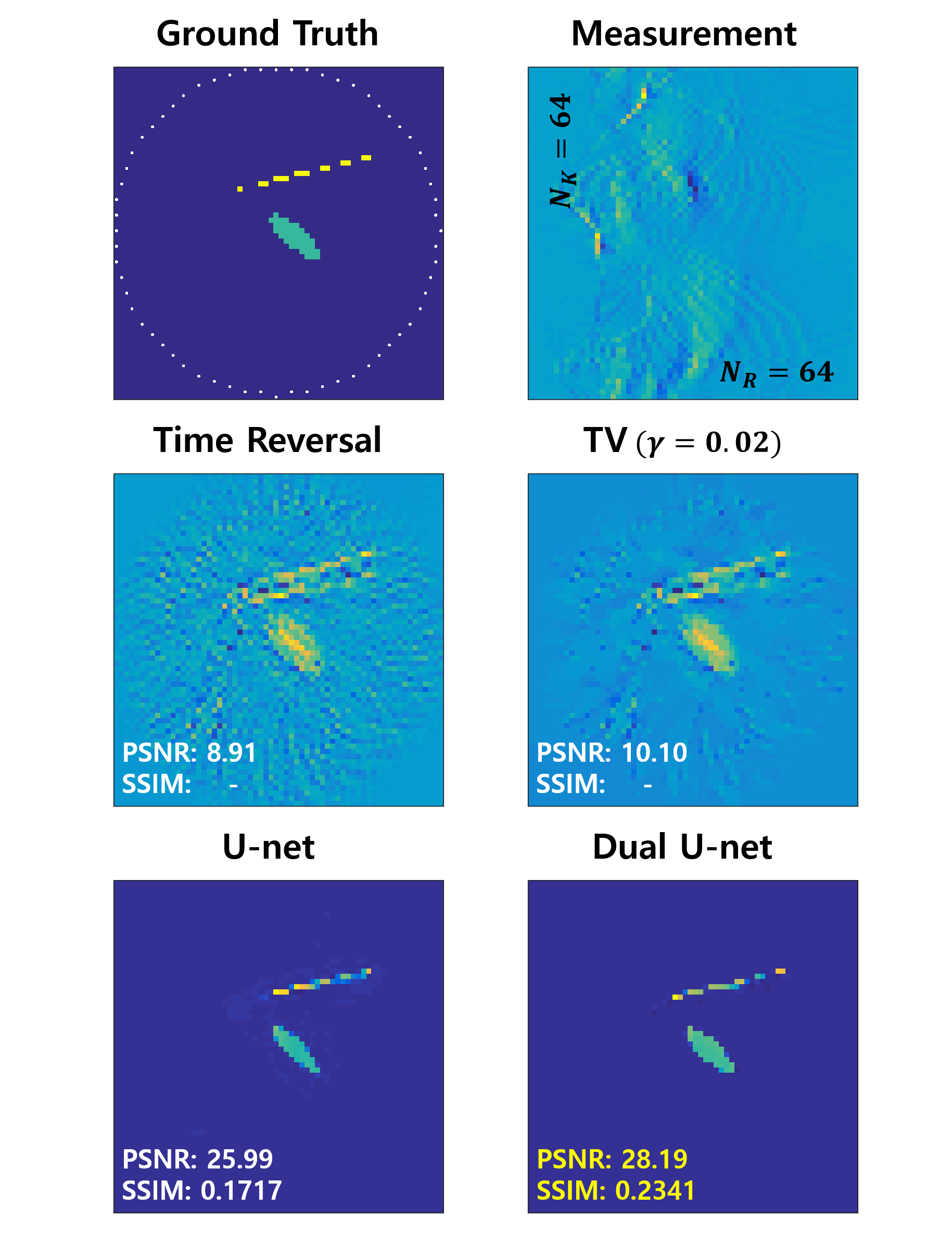}
\caption{The denoised test data-set images using various algorithms in a clean measurement condition. For a fair comparison, we normalized the pixel values to lie in $[0,1]$.}
\label{fig:res-els}
\end{figure}

\begin{figure}[!hbt]
\centering
\includegraphics[width=0.8\linewidth]{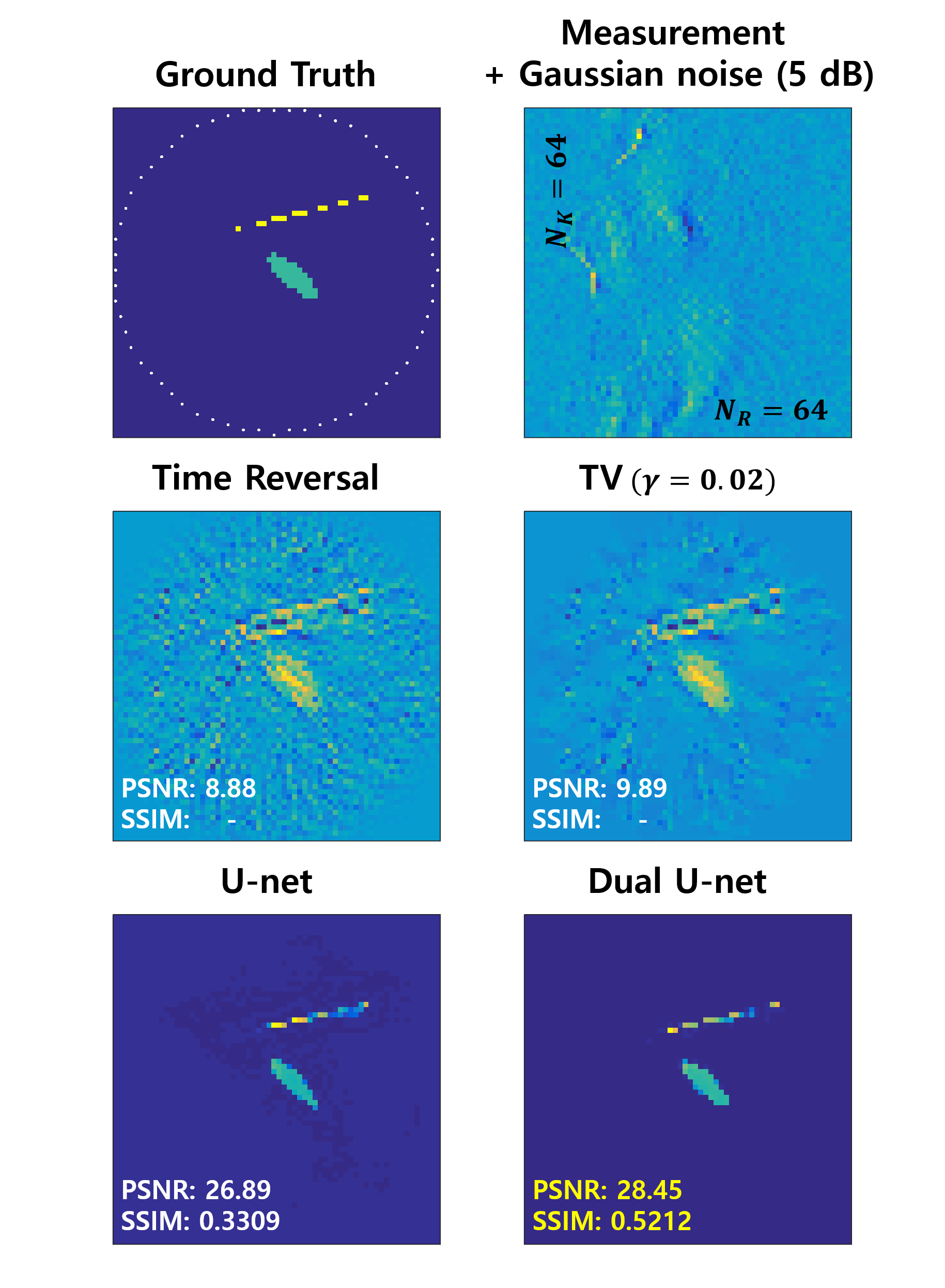}
\caption{The denoised test data-set images using various algorithms in a noisy measurement condition. For a fair comparison, we normalized the pixel values to lie in $[0,1]$.}
\label{fig:res-els-noise}
\end{figure}

\begin{figure}[!hbt]
\centering
\includegraphics[width=0.8\linewidth]{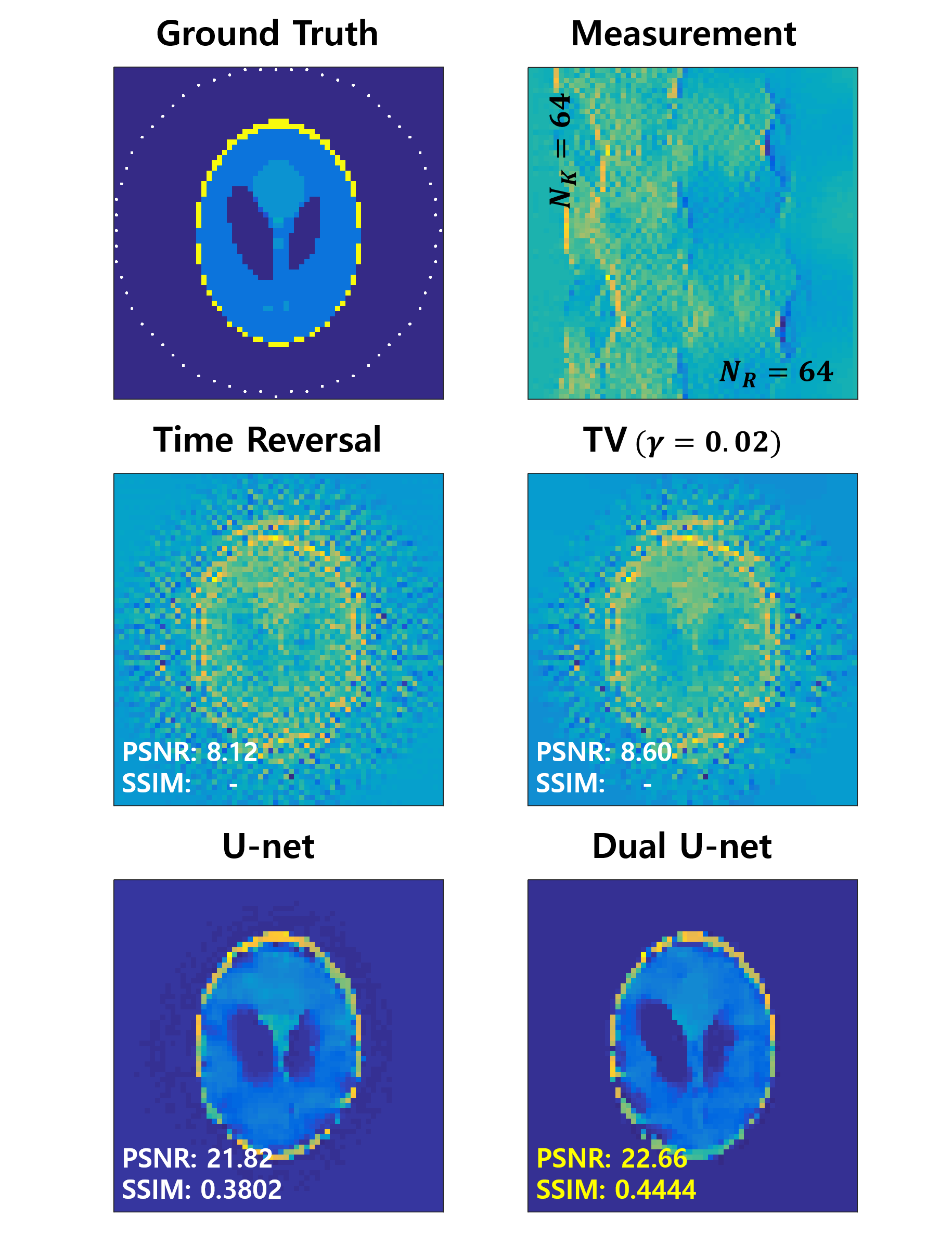}
\caption{The denoised Shepp-Logan phantom images using different algorithms in a clean measurement condition. For a fair comparison, we normalized the pixel values to lie in $[0,1]$.}
\label{fig:res-sl}
\end{figure}

\begin{figure}[!hbt]
\centering
\includegraphics[width=0.8\linewidth]{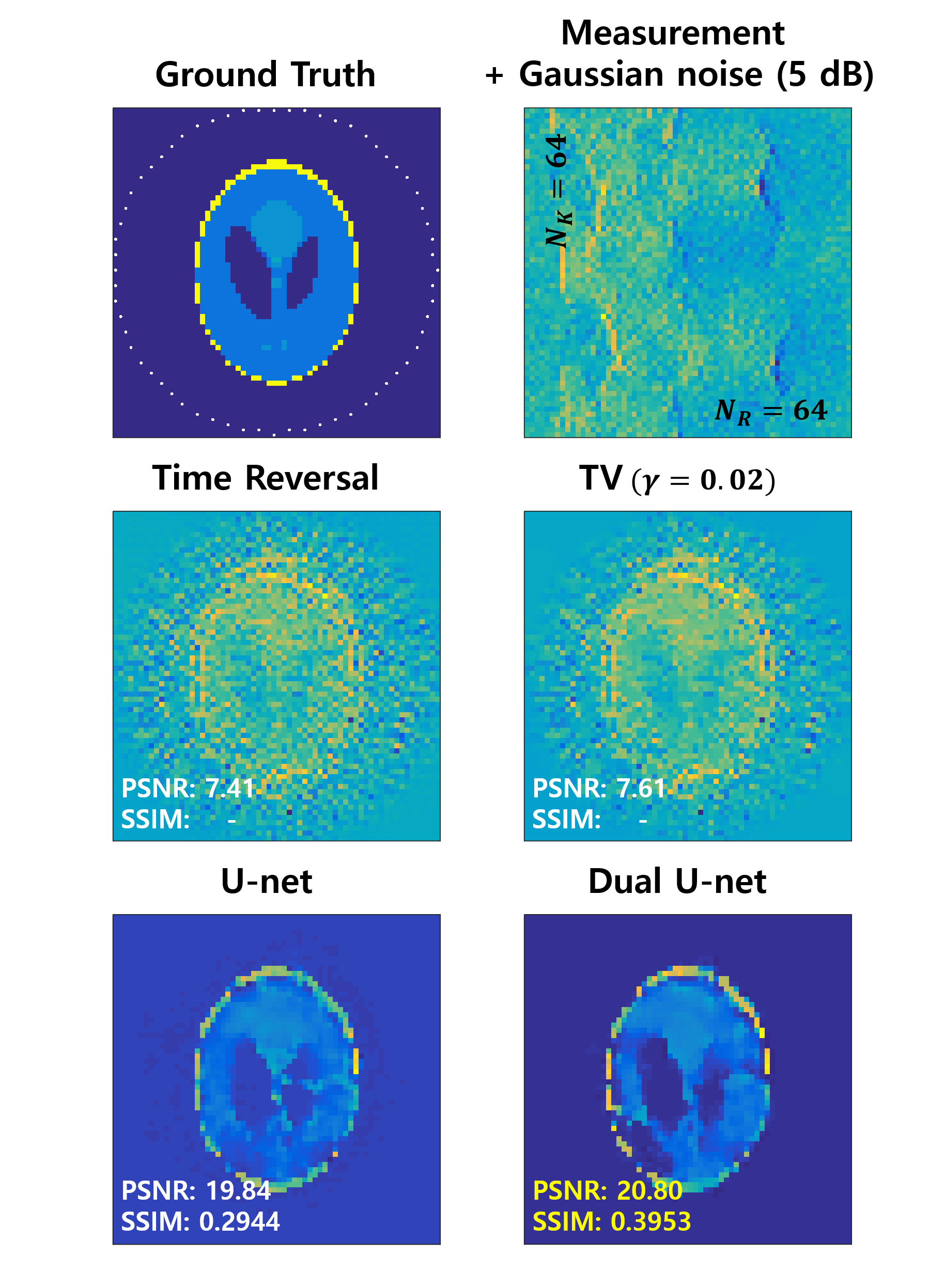}
\caption{The denoised Shepp-Logan phantom images using different algorithms in a noisy  measurement condition. For a fair comparison, we normalized the pixel values to lie in $[0,1]$.}
\label{fig:res-sl-noise}
\end{figure}

\Cref{fig:res-els,fig:res-els-noise}, and \Cref{fig:res-sl,fig:res-sl-noise} show the reconstruction results from the test data-set and Shepp-Logan phantom data.
In particular, \Cref{fig:res-els,fig:res-sl} correspond to the noiseless measurements, while
\Cref{fig:res-els-noise,fig:res-sl-noise}  are from the noisy measurements.
In a noisy condition, a white Gaussian noise with SNR$=5dB$ is added to the measurement and the images are reconstructed using the time-reversal algorithm. 
 In both conditions, for the total variation algorithm, a regularization parameter $\gamma=0.02$ is chosen without any other constraint and  the FISTA algorithm is used. Here, we could not find a significant improvement in the quality of the images by varying the hyperparameter $\gamma$. 
 These results are compared with the results by  the neural networks which are trained on the images from the clean measurements only.  
 Note that the network has seen neither the images from the noisy measurements nor  the Shepp-Logan phantom during the training phase. 

\subsection{Discussion}
\begin{figure}[!hbt]
\centering
\includegraphics[width=0.8\linewidth]{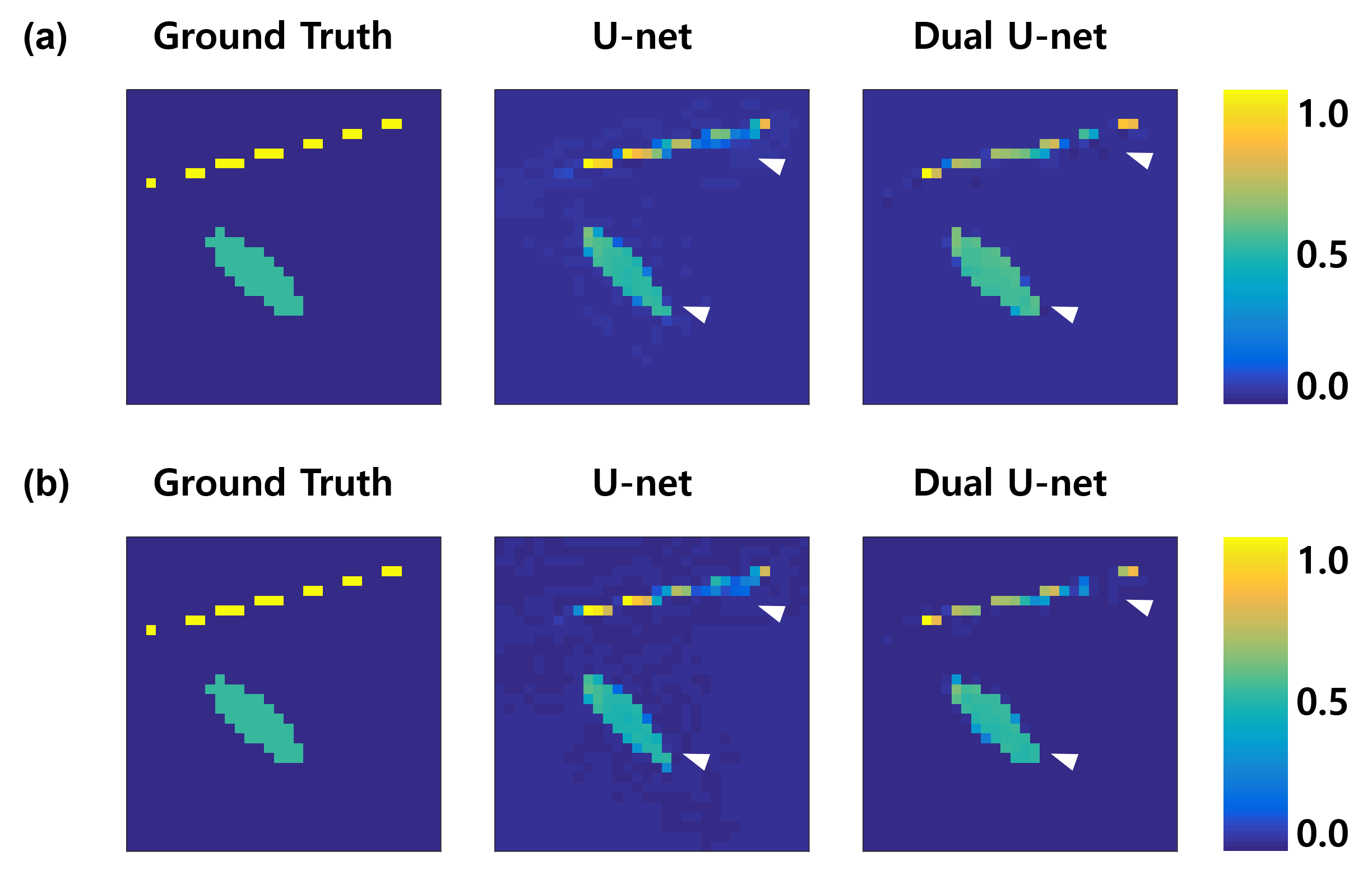}
\caption{The zoomed-in versions of denoised test data-set images in both (a) clean and (b) noisy measurement conditions. 
}
\label{fig:res-els-zoom}
\end{figure}

\begin{figure}[!hbt]
\centering
\includegraphics[width=0.8\linewidth]{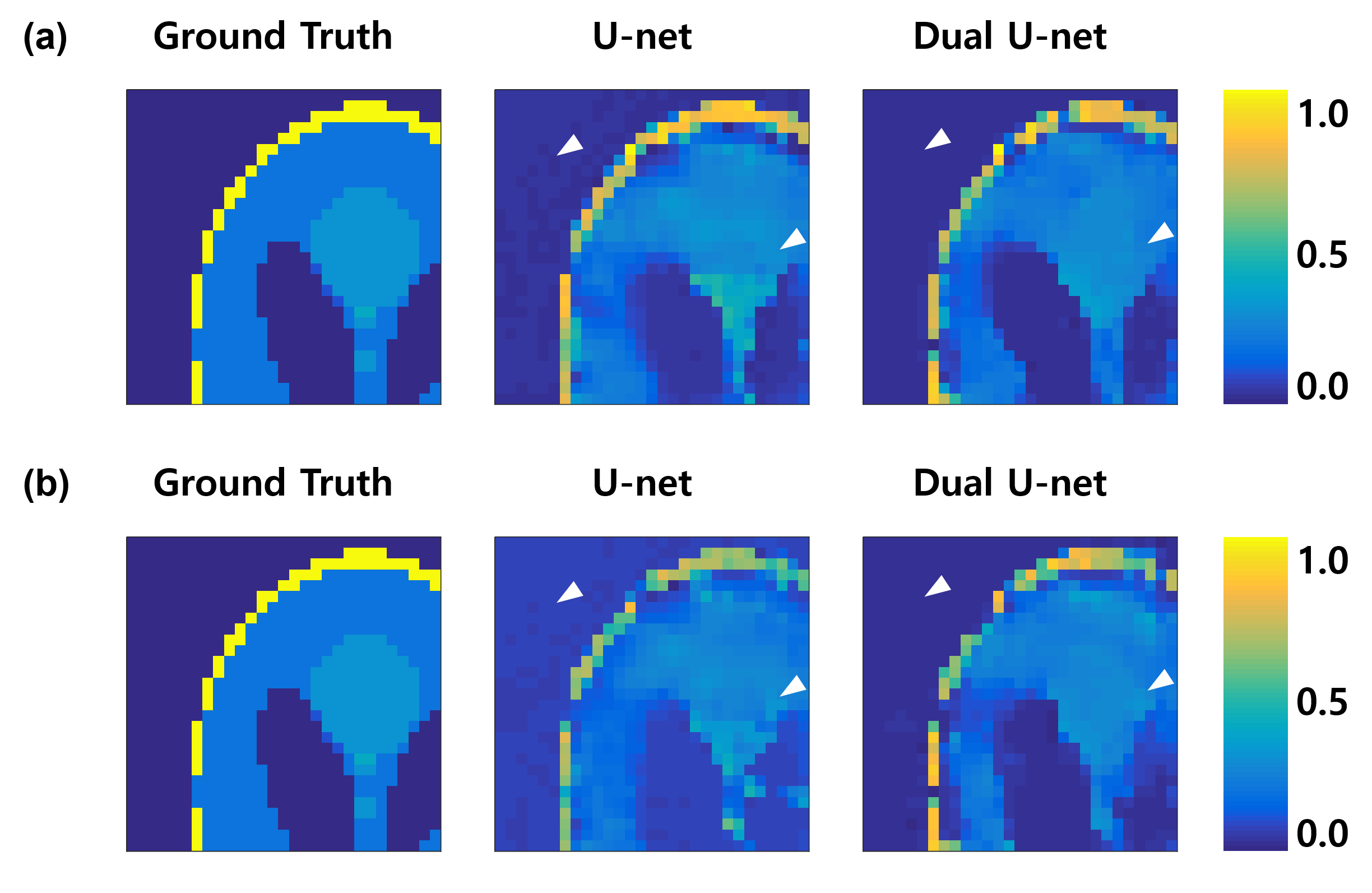}
\caption{The zoomed-in versions of denoised Shepp-Logan phantom images in both (a) clean and (b) noisy measurement conditions. 
}
\label{fig:res-sl-zoom}
\end{figure}

The denoising methods using the neural networks showed a superior performance over the total variation algorithm. Among those, the  Dual-Frame U-Net showed the best results in both PSNR and SSIM. Though the standard U-Net recovered the overall shapes of the inclusions, it failed to find an accurate outfit and lacks the fine details of the inclusions (see \Cref{fig:res-els-zoom,fig:res-sl-zoom}). For example, the recovered shapes of the ellipses using Dual-Frame U-Net in thin and sparse inclusions case have sharper ends than standard U-Net relative to that of the ground truth as highlighted in \Cref{fig:res-els-zoom}. In addition, the standard U-Net failed to remove artifacts around the inclusions and bias in the background (\Cref{fig:res-els-zoom,fig:res-sl-zoom}). On the other hand, the Dual-Frame U-Net recovered the oval shapes of the inclusions and their pixel values more accurately in both sparse and extended targets (pointed out by white arrows in \Cref{fig:res-els-zoom,fig:res-sl-zoom}). These differences come from the overly emphasized low frequency components in the U-Net configuration  that does not meet the frame condition (see \cref{ss:DUN}). 

\section{Conclusion}\label{Sect:Con}
In this article, we showed that the problem of elastic source imaging with very sparse data, both in space and time, can be successfully dealt with our proposed deep learning framework. While the conventional denoising algorithm using TV regularization gives an unsatisfying reconstruction quality, deep learning approaches showed more robust reconstruction with better peak signal-to-noise ratio (PSNR) and structural similarity index (SSIM). We showed that the network performance can be further improved by using the Dual-Frame U-Net architecture, which satisfies a frame condition. 

\section*{Acknowledgments}
The authors would like to thank Dr. Elie Bretin for providing the source code for the forward elastic solver and weighted time-reversal algorithm.

\section*{Appendix}

\begin{figure}[!hbt]
\centering
\includegraphics[width=12cm]{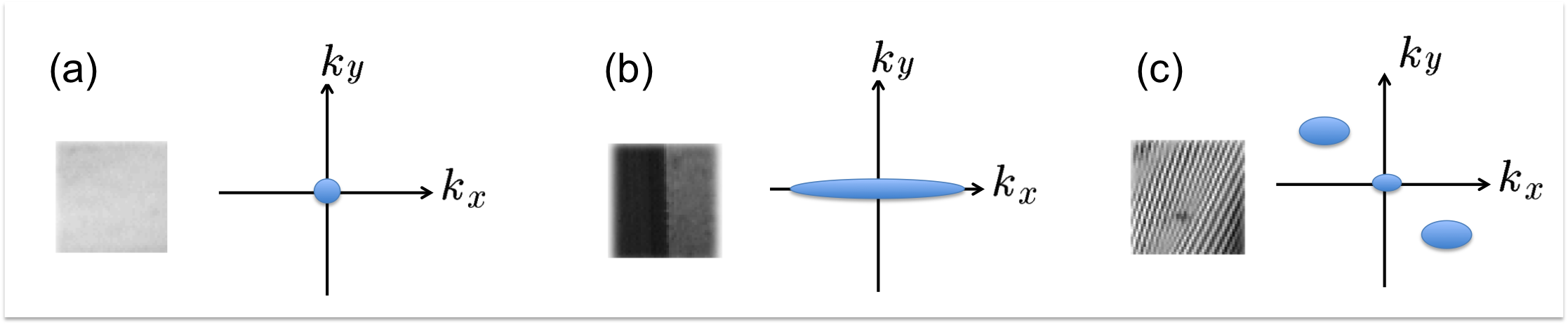}
\caption{Spectral components of patches. (a) Smooth background:  spectral components are mostly concentrated in the low frequency regions.
(b) Edge patch: spectral components are elongated perpendicular to the edge. (c) Texture patch: spectral components are distributed at  the harmonics  of the texture orthogonal to texture direction. }
\label{fig:flowchart}
\end{figure}

To make this paper self-contained,  here we briefly review the origin of the low-rank Hankel matrix as extensively studied
in \cite{ye2017deep,ye2016compressive}.

Note that many types of image patches have sparsely distributed Fourier spectra.
For example,   as shown in \Cref{fig:flowchart}(a),  a smoothly varying patch usually has spectrum content in the low-frequency regions.
For the case of an edge as shown in \Cref{fig:flowchart}(b), the spectral components are mostly localized along the $k_x$-axis. 
Similar spectral domain sparsity can be observed in the texture patch shown in  \Cref{fig:flowchart}(c), where the spectral components of the patch are distributed at the harmonics of the texture.
In these cases,  if we construct a Hankel matrix using the corresponding image patch,  the resulting Hankel matrix is low-ranked \cite{ye2016compressive}.

In order to understand this intriguing relationship,  consider a 1-D signal, whose spectrum  in the Fourier domain is sparse and can
be modeled as the sum of Dirac masses:
\begin{equation}\label{eq:signal3}
\hat f(\omega) = 2\pi \sum_{j=0}^{r-1} c_{j} \delta \left( \omega- \omega_j \right), \qquad \omega_j \in [0, 2\pi] ,
\end{equation}
where $\{\omega_j\}_{j=0}^{r-1}$  refers to the corresponding sequence of  the harmonic components in the Fourier domain.
Then, the corresponding discrete time-domain signal is  given by:
 \begin{eqnarray}\label{eq:fs}
[\fb]_k = \sum_{j=0}^{r-1} c_{j} e^{-i k \omega_j }.
\end{eqnarray}
Suppose that  we have a $(r+1)$-length  filter $\mathbf{h}$ that has the  $z$-transform representation  \cite{vetterli2002sampling}  
\begin{eqnarray}\label{eq:afilter}
\hat h(z)  &=& \sum_{l=0}^r  [\mathbf{h}]_l z^{-l} = \prod_{j=0}^{r-1} (1- e^{-i\omega_j} z^{-1}) \ .
\end{eqnarray}
Then, it is easy to see that \cite{vetterli2002sampling}
\begin{eqnarray}\label{eq:annf}
 \fb\circledast \mathbf{h}  =\mathbf{0}.
\end{eqnarray}
Thus, the filter $\mathbf{h}$ annihilates the signal $\fb$ and is accordingly referred to as the \emph{annihilating filter}.
Moreover,  since Eq.~\eqref{eq:annf} can be represented as
$$
\hank_p(\fb) \mathbf{h}' =\mathbf{0},
$$
the Hankel matrix $\hank_p(\fb)$ is rank-deficient.   In fact, the rank of the Hankel matrix can be explicitly determined by the size of the minimum-size annihilating filter \cite{ye2016compressive}. Therefore, if the matrix pencil size $p$ is chosen bigger than the minimum annihilating filter size, the Hankel matrix is low-ranked.

\bibliographystyle{siamplain}

\end{document}